\documentclass[amssymb,12pt]{article}

\usepackage{amsmath}
\usepackage{graphicx,
}  
\usepackage[psamsfonts]{amssymb}
\usepackage{amscd}

\title{$C^1$ genericity of trivial centralizers}
\author{Christian Bonatti, Sylvain Crovisier and Amie Wilkinson}

 \def\NN{{\mathbb N}}  
 \def\RR{{\mathbb R}}  \def\TT{{\mathbb T}}
   
 \def\ZZ{{\mathbb Z}}

\def\La{\Lambda}
\def\De{\Delta}
\def\Om{\Omega}
\def\Ga{\Gamma}

\def\cC{{\cal C}}  \def\cI{{\cal I}} \def\cO{{\cal O}} \def\cU{{\cal U}}
\def\cD{{\cal D}}   \def\cP{{\cal P}} \def\cV{{\cal V}}
\def\cE{{\cal E}}    
  \def\cL{{\cal L}} \def\cR{{\cal R}}

\newtheorem{theorem}{Theorem}[section]
\newtheorem{proposition}[theorem]{Proposition}
\newtheorem{lemma}[theorem]{Lemma}
\newtheorem{corollary}[theorem]{Corollary}

\newtheorem{theo}[theorem]{Theorem}
\newtheorem{clai}{Claim}
\newtheorem{lemm}[theorem]{Lemma}
\newtheorem{coro}[theorem]{Corollary}
\newtheorem{defi}[theorem]{Definition}
\newtheorem{prop}[theorem]{Proposition}

\newtheorem{rema}[theorem]{Remark}
\newtheorem{ques}[theorem]{Question}

\newtheorem{affi}{Claim}

\newenvironment{demo}[1][{\noindent Proof}]{{\bf #1. }}{\hfill$\Box$\medskip}

\title{Centralizers of $C^1$-generic diffeomorphisms}
\author{C. Bonatti, S. Crovisier and A. Wilkinson}

\def\sign{\hbox{Sign} }

\def\endproof{\hfill$\Box$\medskip}

\def\Diff{\hbox{Diff} }

\def\eproof{\hfill$\Box$ \bigskip}

\def\title{\em}
\def\eps{\varepsilon}

\def\bar{\overline}

\def\Jac{\hbox{Jac}}

\def\T{\mathcal{T}}

\def\transverse{\,\raise2pt\hbox to1em{\hfil$\top$\hfil}\hskip -1em \hbox
to1em{\hfil$\cap$\hfil}\,} 
\newcommand\R{\mbox{\bf R}}

\newlength{\figboxwidth} 

\begin{document}

\maketitle
\begin{abstract}
On the one hand, we prove that the spaces of $C^1$ symplectomorphisms
and of $C^1$ volume-preserving diffeomorphisms 
both contain residual subsets of diffeomorphisms
whose centralizers are trivial. On the other hand, we show that
the space of $C^1$ diffeomorphisms of the circle and a non-empty open set of $C^1$ diffeomorphisms
of the two-sphere contain dense subsets of diffeomorphisms
whose centralizer has a sub-group isomorphic to $\RR$.

\vskip 2mm

\begin{description}
\item[\bf Key words:] Trivial centralizer, trivial symmetries, 
Mather invariant, $C^1$ generic properties.
\end{description}
\end{abstract}

\section*{Introduction}

Let $M$ be a connected compact manifold.
The {\em centralizer} of a $C^r$ diffeomorphism $f\in \Diff^r(M)$ is defined 
as $$C(f):=\{g\in \Diff^r(M): fg=gf\}.$$
Clearly $C(f)$ always contains the group $<f>$ of all the powers of $f$.  We 
say that $f$ has {\em trivial centralizer} if $C(f) = <f>$.  A
diffeomorphism $f$ with trivial centralizer posesses no smooth symmetries,
such as those that would arise if, for example, $f$ embedded in a flow or
were the lift of another diffeomorphism. Smale asked the following:

\begin{ques}[\cite{Sm1,Sm2}]\label{c=main} Let 
$\T^r(M)\subset \Diff^r(M), r\geq 1$ denote the set of
$C^r$ diffeomorphisms of a compact manifold $M$ with trivial centralizer.  

\begin{enumerate} 
\item Is $\T^r(M)$ dense in $\Diff^r(M)$?
\item Is $\T^r(M)$ residual in $\Diff^r(M)$? 
\item Is $\T^r(M)$ open in $\Diff^r(M)$?
\end{enumerate}
\end{ques}

This question has been answered in several special cases.  To
summarize these results in rough chronological order, we have:
\begin{itemize} 
\item $\T^r(S^1)$ is open and dense in $\Diff^r(S^1)$ for
  $r\geq 2$ \cite{Ko};
\item  $\T^1(M)$ is residual among the Axiom A diffeomorphisms in $\Diff^1(M)$;
in particular, $\T^1(S^1)$ is residual in $\Diff^1(S^1)$ \cite{To1, To2};
\item $\T^\infty(M)$ is open and dense among the Axiom A
  diffeomorphisms in $\Diff^\infty(M)$ 
possessing at least one periodic sink or source \cite{PY1};
\item  $\T^\infty(M)$ is open and dense among the Anosov
 diffeomorphisms in $\Diff^\infty (\TT^n)$, where $\TT^n$ is the
 $n$-torus 
\cite{PY2}; 
\item $\T^\infty(M)$ is locally residual among the partially
hyperbolic diffeomorphisms with $1$-dimensional center \cite{Bu}.
\end{itemize}

There are two main results in this paper.
In the first (Theorem~\ref{t=main1}), 
we give a complete answer to  the first two parts of
Question \ref{c=main} for all compact $M$ in the case of volume-preserving 
and symplectic
$C^1$-diffeomorphisms.  In the second result (Theorem~\ref{t=main2}),
we answer the third part of Question \ref{c=main} for the circle $S^1$
and the sphere $S^2$, again in the case $r=1$.

\subsection*{A) Trivial centralizer for 
$C^1$-generic symplectomorphisms and volume-preserving diffeomorphisms}

In order to state our first main result precisely we will 
need some notation.  If $M$ carries a volume $\mu$, then we denote by $\Diff^1_{\mu}(M)$
the space of $C^1$ diffeomorphisms of $M$ that preserve $\mu$.  If $M$
is a symplectic manifold, then $\hbox{Symp}^1(M)$ denotes the space of
$C^1$ symplectomorphisms of $M$.  The spaces $\Diff^1(M)$,
$\Diff^1_{\mu}(M)$, $\hbox{Symp}^1(M)$ are Baire spaces in the $C^1$
topology.  Recall that a {\em residual} subset of a Baire space is one
that contains a countable intersection of open-dense sets.

\begin{theo}\label{t=main1} Let $M$ be a compact, connected manifold
  of dimension at least $2$.  Then:
\begin{enumerate}
\item[(a)] $\T^1(M)\cap \Diff^1_{\mu}(M)$ is residual in $\Diff^1_{\mu}(M)$.
\item[(b)] $\T^1(M)\cap\hbox{Symp}^1(M)$ 
is residual in ${Symp}^1(M)$.
\end{enumerate}
\end{theo}
Theorem~\ref{t=main1} is a corollary of parts (b) and (c) of the following result:
\begin{theo}\label{t=main} (a) There is a residual set $\cR\subset \Diff^1(M)$ such that,
for any diffeomorphism $f\in \cR$, for any $g\in C(f)$ and any periodic point $x\in Per(f)$, the point $x$ is hyperbolic and,  there exist $m, n\in\ZZ$ such that $g$ coincides with $f^n$ on $W^s(x)$ and with $f^m$ on $W^u(x)$.

\noindent (b) There is a residual set $\cR_{\hbox{symp}}\subset \hbox{Symp}^1(M)$ such that, for any diffeomorphism $f\in \cR_{\hbox{symp}}$ for any $g\in C(f)$ and any hyperbolic periodic point $x\in Per(f)$,   there exist $m, n\in\ZZ$ such that $g$ coincides with $f^n$ on $W^s(x)$ and with $f^m$ on $W^u(x)$.

\noindent(c) There is a residual set $\cR_\mu\subset \Diff^1_{\mu}(M)$ such that, for any diffeomorphism $f\in \cR_\mu$ for any $g\in C(f)$ and any hyperbolic periodic point $x\in Per(f)$,  there exists $n\in\ZZ$ such that $g$ coincides with $f^n$ on 
either $W^s(x)$ or $W^u(x)$.

\end{theo}
Theorem~\ref{t=main}~(a) was previously proved by Togawa \cite{To1,
  To2}, using different methods.    Togawa's methods, combined with
  the results in Appendix A, can also be
  used to prove parts (b) and (c) of Theorem~\ref{t=main}.  While
  using Togawa's results would shorten considerably the proof of
  Theorem~\ref{t=main}, we believe our approach, in particular
  Propositions~\ref{p=main} and \ref{p=diskdist}, has independent
  interest. It would be interesting to see if these results have
  further application.  We discuss the motivation and background to this approach
  in Section~\ref{s=distortion}.

\bigskip

\noindent{\bf Proof of Theorem~\ref{t=main1}.}
Theorem~\ref{t=main1} follows immediately from Theorem~\ref{t=main} and:
\begin{theo}[\cite{BC,ABC}] For any compact connected manifold $M$, 
there are residual sets 
$\widetilde\cR_\mu\subset \Diff^1_\mu(M)$ and $\widetilde\cR_{\hbox{symp}}\subset \hbox{Symp}^1(M)$ such that,  every $f\in \widetilde\cR_\mu\cup \widetilde\cR_{\hbox{symp}}$ has a hyperbolic periodic point $p$ with
$$\overline{W^s(p)} = \overline{W^u(p)} = M.$$ 
\end{theo}\eproof

More generally, Theorem~\ref{t=main} naturally applies to the class of 
$C^1$ diffeomorphisms satisfying a property we call periodic 
accessibility.
A diffeomorphism $f$ on a compact manifold satisfies the  
{\em periodic accessibility property}
if there is a dense subset $\cE\subset M$ of non-periodic points such that any pair of points $x,y\in \cE$ may be joined by a finite sequence
$x_0=x,x_1,\dots,x_n=y$, $x_i\in \cE$ and a sequence $p_i$ of hyperbolic periodic orbits such that
 for any $i\in\{0,\dots,n-1\}$ one has:
$$\{x_i,x_{i+1}\} \subset \overline{W^s_{orb}(p_i)} \quad \mbox{ or }\quad  \{x_i,x_{i+1}\} \subset \overline{W^s_{orb}(p_i)}.$$ 

\begin{ques} Is periodic accessibility generic in $\Diff^r(M)$?
\end{ques}
As a weaker problem, one can also ask if, for generic diffeomorphims,
the union of the stable manifolds of the periodic points are dense in $M$.
\bigskip

Theorem~\ref{t=main} has the immediate corollary:
\begin{coro} 
Furthermore, if $f\in \cR$ satisfies the periodic accessibility property then 
$C(f)$ is trivial.
\end{coro}

The periodic accessibility property is satisfied by Axiom A
diffeomorphisms, by $C^1$-generic tame diffeomorphisms (i.e. by
$C^1$-generic diffeomorphisms having finitely many homoclinic
classes), and by $C^1$-generic conservative (volume preserving or
symplectic) diffeomorphisms. In this way, one can recover Togawa's
result that the $C^1$-generic Axiom A diffeomorphism has trivial centralizer.

\subsection*{B) Large centralizer for a locally $C^1$ dense set of diffeomorphisms}

Our next main result addresses the third part of Question~\ref{c=main}: is
$\T^r(M)$ open in $\Diff^r(M)$?   In the case of the circle, recall that Kopell
proved that $\T^r(S^1)$ is open-dense in $\Diff^r(S^1)$ for $r\geq 2$,
and Togawa proved that $\T^1(S^1)$ is residual in $\Diff^1(S^1)$.  
It is natural to ask whether Togawa's result can be strengthened to
show that $\T^1(S^1)$ is open-dense.

Our next main result
shows that the answer is ``no'': the answers to the third part of Question~\ref{c=main}
are genuinely different in  
the $C^1$ and $C^2$ topologies, at least for the circle.  We
are also able to answer the third part of  Question~\ref{c=main} for
the $2$-sphere. Specifically, we have:

\begin{theo}\label{t=main2}  $\T^1(S^1)$ and  $\T^1(S^2)$ are not open. Moreover:
 
\medskip

\noindent (a) There is a dense subset $\cD^1\subset\Diff^1(S^1)$ such that every $f\in\cD^1$  leaves invariant a $C^\infty$ Morse-Smale vector field. 
In particular, $C(f)$ contains a subgroup isomorphic to $\RR$. 

\medskip

\noindent (b) Let $\cO\subset \Diff^1(S^2)$ denote the (open) 
subset of Morse-Smale diffeomorphisms 
$g$ such that the nonwandering set 
$\Om(g)$ consists of two fixed points, one source $N_g$ and one sink $S_g$, 
such that the derivatives $D_{N_g}g$ and $D_{S_g}g$ have each a complex 
(non real) eigenvalue. 

There is a dense subset $\cD^2\subset\cO$ such that every $f\in\cD^2$ is the time  $1$ map of a Morse-Smale $C^\infty$-vector field. In particular, $C(f)$ contains a subgroup  isomorphic to $\RR$. 

\end{theo}

\subsection*{Structure of the paper}
In order to prove Theorem~\ref{t=main}, it is enough to show that along the invariant manifolds of the periodic
points, $f$ satisfies an unbounded distortion property. This is discussed in Section~\ref{s=distortion}.
As a simpler setting, we also deal with contractions of $\RR^d$ whose unique periodic point is $0$.
In Section~\ref{s=proof-main}, we will see that $C^1$-generic contraction of $\RR^d$
has the unbounded distortion property; this can be generalized to the dynamics inside the invariant manifolds
of the periodic points, since by Appendix~\ref{a=invmanif}, any perturbation of the dynamics inside the
stable manifod of a periodic point can be realized as a perturbation of the dynamics on $M$.
Theorem~\ref{t=main2} will be proved in Section~\ref{s=contre-exemple}.

\section{The unbounded distortion property}\label{s=distortion}
Kopell's proof in \cite{Ko} that $\T^r(S^1)$ is open-dense in $\Diff^r(S^1)$ for
$r\geq 2$  uses the fact that
a $C^2$ diffeomorphism $f$ of $[0,1]$  without fixed points in $(0,1)$
has {\em bounded distortion}, meaning: for any $x,y\in (0,1)$, the ratio
\begin{eqnarray}\label{e=xydistort}
\frac{|{f^n}'(x)|}{|{f^n}'(y)|}
\end{eqnarray}
is bounded, independent of $n$ and uniformly for $x,y$ lying in a
compact set.  A bounded distortion estimate 
lies behind many results about $C^2$, hyperbolic 
diffeomorphisms of the circle and 
codimension-$1$ foliations.

Suppose that $r\geq 2$. 
Since Morse-Smale diffeomorphisms are open and dense in $\Diff^r(S^1)$,
the proof that $\T^r(S^1)$ is open-dense in $\Diff^r(S^1)$ 
essentially reduces to showing that ($C^r$-open and densely)
a $C^r$ diffeomorphism $f:[0,1]\to [0,1]$  without fixed points in $(0,1)$
has trivial centralizer.  The 
bounded distortion of such an $f$ forces its centralizer to
embed simultaneously in two smooth flows containing $f$, one determined
by the germ of $f$ at $0$, and the other by the germ at $1$; 
for an open and dense set of $f\in \Diff_+^r[0,1]$, these flows
agree only at the iterates of $f$.
The $r\geq 2$ hypothesis is clearly necessary for bounded distortion.

The central observation and starting point of this paper is that
the centralizer of a $C^1$ diffeomorphism of $[0,1]$ with 
{\em unbounded} distortion is always trivial.  We elaborate a bit on this.
Notice that if $x$ and $y$ lie on the same $f$-orbit, then the 
ratio in (\ref{e=xydistort}) is bounded, independent of $n$. 
We show that, $C^1$-generically among the diffeomorphisms of $[0,1]$ without
fixed points $(0,1)$, the ratio (\ref{e=xydistort}) is uniformly bounded in $n$ {\em only if}
$x$ and $y$ lie on the same orbit; that is, for a residual set of $f$,
and for all $x,y\in (0,1)$, 
if $x\notin {\mathcal O}_f(y) = \{f^n(y)\,\vert\, n\in \ZZ\}$,
then
\begin{eqnarray}\label{e=xylimsup}
\limsup_{n\to\infty} \frac{|{f^n}'(x)|}{|{f^n}'(y)|} = \infty.
\end{eqnarray}

 Assume that this unbounded  distortion property holds for $f$.  
Fix $x\in(0,1)$.  A simple application of the Chain
Rule shows that if $gf = fg$, then the distortion in (\ref{e=xydistort}) 
between  $x$ and $y=g(x)$ is bounded; hence $x$ and $g(x)$ must
lie on the same $f$-orbit. From here, it is straightforward to show that 
$g= f^n$, for some $n$ (see Lemma~\ref{l=fixespoint} and
Corollary~\ref{c=fixesorbit} below).  
As in \cite{Ko}, a small
amount of additional work shows that a residual set in $\Diff^1(S^1)$
has trivial centralizer.  The details of this argument we have just 
described for $S^1$ are contained in this section and 
Section~\ref{s=linconf}.   

The bulk of this paper is devoted to formulating and 
proving a higher-dimensional version of the argument we have just 
described. The interval is replaced by an invariant manifold (stable or 
unstable) of a periodic point.  The derivative $f'$ in 
(\ref{e=xylimsup}) is replaced by
the Jacobian of $f$ along the invariant manifold.

\subsection{Unbounded distortion along invariant manifolds}

Let $f:M\to M$ be a $C^1$ diffeomorphism, and let $p\in M$ be a hyperbolic
periodic point of $f$.  For $x\in W^s(p)$ we denote by $\Jac^s(f)(x)$
the Jacobian of the map induced by $T_xf$ between $T_xW^s(p)$ and
$T_{f(x)}W^s(f(p))$. 

\begin{defi} A hyperbolic periodic point $p\in M$
has the {\em stable manifold  distortion property} if, for every
$x,y\in W^s(p)\setminus \{p\}$ not in the same $f$-orbit,
$$\limsup_{n\to\infty} \left|\frac{\Jac^s(f^n)(x)}{\Jac^s(f^n)(y)}\right| = \infty.$$
\end{defi}

As mentioned in the previous subsection, unbounded distortion forces trivial 
centralizers:

\begin{lemma}\label{l=fixespoint} Let $p$ be a hyperbolic periodic point of period $k\in \NN$
with the stable manifold distortion property,
and let $g\in C(f)$.  If $g(p) = p$, then there exists an $m\in \ZZ $
such that $g = f^{km}$ on $W^s(p)$.
\end{lemma}
\noindent{\bf Proof of Lemma~\ref{l=fixespoint}.} We claim that for every 
$x\in W^s(p)$, $g$ preserves the $f^k$-orbit of $x$.  From this claim it
follows that for every $x\in W^s(p)$, there exists an integer $m(x)$ 
such that $g(x) = f^{km(x)}(x)$, and there is
a unique such $m(x)$ if $x\neq p$.  Continuity of $f$ implies that the function
$m$ is locally constant on $W^s(p)\setminus \{p\}$.   If $\dim(W^s(p))>1$, then $W^s(p)\setminus \{p\}$ is connected, and
$m$ is constant.  If $\dim(W^s(p))=1$, then $m$ is constant on each of the
connected components of  $W^s(p)\setminus \{p\}$; in this case,
since $g$ is differentiable
at $p$ and $|{f^k}'(p)|\neq 1$, the values of $m$ on the two components must coincide.

It remains to prove the claim. We may assume that $x\neq p$. The relation
$gf^n = f^n g$ implies that,
$$\Jac^s(g)(f^n x)\Jac^s(f^n)(x) =  \Jac^s(f^n)(g x)\Jac^s(g)(x).$$
In particular, for all $m\geq 0$, we have
$$\left|\frac{\Jac^s(g)(f^{mk} x)}{\Jac^s(g)(x)}\right| = 
\left|\frac{\Jac^s(f^{mk})(g x)}{\Jac^s(f^{mk})(x)}\right|.$$
Since $f^{mk}(x)$ lies in a compact region of $W^s(p)$ for all $m\in \NN$,
the left hand side of this expression is uniformly bounded in $m$.  On the other hand, it is easy to see that if $f$ has the stable manifold distortion property, then so does $f^k$.  This implies that the right hand side of the equation
above is unbounded, a contradiction.  This proves the claim.\eproof

\begin{corollary}\label{c=fixesorbit} Let $p$ be a hyperbolic periodic point of period $k\in \NN$
with the stable manifold distortion property,
and let $g\in C(f)$.  If $g({\mathcal O}(p,f)) = {\mathcal O}(p,f)$, then there exists an $m\in \ZZ $
such that $g = f^{m}$ on $W^s({\mathcal O}(p,f))$.
\end{corollary}

\noindent{\bf Proof of Corollary~\ref{c=fixesorbit}.} Since $g$ preserves the $f$-orbit of $p$, we have $g(p) = f^j(p)$ for some integer $j$.  Let $G = f^{-j}g$. Then it is easy to see that $G$ commutes with $f$, and $G$ fixes every point 
on the $f$-orbit of $p$.  By Lemma~\ref{l=fixespoint}, there is an integer $m$ such that the restriction of
$G$ to $W^s(p)$ coincides with $f^{mk}$.  Again, since $f$ and
$G$ commute, the restriction of $G$ to
$W^s(f^ip)$ is conjugate by $f^i$ to the restriction of $G$ to 
 $W^s(p)$.  Consequently, $G$ coincides with $f^{mk}$ on the stable manifold 
 $W^s({\mathcal O}(p,f))$ of the orbit of $p$, and so $g$ coincides
with $f^{mk+j}$ on  $W^s({\mathcal O}(p,f))$.\eproof

Note that for any integer $k\geq1$,
a $C^1$-generic diffeomorphism $f$ has only finitely periodic orbits of period $k$
and, by transversality, that all these orbits have different exponents. In particular, any
diffeomorphism $g\in C(f)$ preserves each of these orbits and satisfies the assumption
of corollary~\ref{c=fixesorbit}.

To prove Theorem~\ref{t=main}, we are thus reduced to proving:

\begin{prop}\label{p=main} There is a residual set $\cR\subset \Diff^1(M)$ such that, for any diffeomorphism $f\in \cR$, every periodic point $x\in Per(f)$ is hyperbolic and has the stable manifold distortion property.

There is a residual set $\cR_{symp}\subset \hbox{Symp}^1(M)$ such that, for any diffeomorphism $f\in \cR_{symp}$, every hyperbolic periodic point $x\in Per(f)$ has the stable distortion property. 

There is a residual set $\cR_{\mu}\subset \Diff^1_\mu(M)$ such that, for any diffeomorphism $f\in \cR_\mu$, and  any hyperbolic periodic point $x\in Per(f)$,
if $W^s(x)$ has codimension at least $\dim(M)/2$, then $x$ has the stable manifold
distortion property.
\end{prop}

\subsection{Contractions of $\RR^d$}
Let $B^d$ denotes the unit closed ball $\overline{B(0,1)}$ of $\RR^d$
and consider the Banach space of $C^1$ maps $B^d\to\RR^d$
that send $0$ to $0$, endowed with the $C^1$-topology given by the $C^1$-norm:
$$\|f-g\|_1=\sup_{x\in B^d} \|f(x)-g(x)\|+ \|D_xf-D_xg\|.$$
The set of embeddings $B^d\to \RR^d$ fixing $0$
defines an open subset that will be denoted by $\cD^d$.
\begin{rema}
Since the origin is  fixed, the metric $\|\cdot\|_1$ is equivalent to the metric
defined by
$$\|f-g\|_1'=\sup_{x\in B^d} \|D_xf-D_xg\|.$$
In general, we will prefer to work with this second one.
\end{rema}

A \emph{contraction} of $\RR^d$ is an element of $\cD^d$ that sends $B^d$ into ${B(0,1)}$,
so that $0$ is a (hyperbolic) sink that attracts all the points in $B^d$.
The set of contractions of $\RR^d$ is an open subset   $\cC^d\subset
\cD^d$, hence a Baire space.

Let $f$ be a diffeomorphism of a manifold $M$, $p$ be a periodic point of $f$ and $n^s$
its stable dimension. A \emph{stable chart} for $p$ is a local chart
$\psi\colon \RR^d\to M$ such that
if one denotes by $\pi$ the projection of $\RR^d$ onto the $n^s$ first coordinates
we have the following properties.
\begin{itemize}
\item The domain $\psi(\RR^d)$ contains $p$.
\item In the chart $\psi$, the local stable manifold of $p$ contains the graph
of a $C^1$ map $g\colon \RR^{n^s}\to \RR^{d-n^s}$.
\item Let $v$ be equal to $\pi(\psi^{-1}(p))$ and let $\theta$ be the $C^1$-map
defined on a neighborhood of $0$ by projecting on the space $\RR^{n^s}$
the dynamics of $f$ in the local stable manifold of $p$:
$$\theta\colon x\mapsto \pi\circ \psi^{-1}\circ f\circ \psi(x+v, g(x+v)) -v,$$
then, $\theta$ belongs to $\cC^{n^s}$.
\end{itemize}

\begin{proposition}\label{p.stable-chart}
Any hyperbolic periodic point $p$ of a diffeomorphism $f$ has a stable chart $\psi$.
Moreover, for any diffeomorphism $g$ in a $C^1$-neighborhood $\cU$ of $f$, the continuation
$p_g$ of $p$ also admits the chart $\psi$ as a stable chart.

The family of contractions $\theta_g$ associated to the periodic point and to the chart $\psi$
induces a continuous map $\Theta\colon \cU\to \cC^{n^s}$. This map is open.
\end{proposition}

In the conservative setting, the same property holds.
\begin{theorem}
Let $\Theta\colon \cU\to \cC^{n^s}$ be a family of contractions associated to
a periodic point $p$ and a stable chart $\psi$ as in proposition~\ref{p.stable-chart}.
Then, the map $\Theta\colon \cU\cap \hbox{Symp}^1(M)\to \cC^{n^s}$ is open.

If the dimension $n^s$ of the stable space of $p$ is larger or equal to $\dim(M)/2$,
then, the map $\Theta\colon \cU\cap \Diff_\mu^1(M)\to \cC^{n^s}$ is open.
\end{theorem}
This will be proved in Sections~\ref{s=symplectic} and \ref{s=volpres}.

\begin{proposition}
For any integer $n\geq 0$, there exists
\begin{itemize}
\item a family $\cP_n$ of pairwise disjoint open subsets whose union is dense in $\Diff^1(M)$,
\item for each $\cU\in \cP_n$, finitely many charts $\psi_1,\dots,\psi_s\colon \RR^d\to M$,
\end{itemize}
such that any diffeomorphism $f\in \cU$ has the following properties:
\begin{itemize}
\item $f$ has $s$ periodic points of period less than $n$, all are hyperbolic.
Each domain $\psi_i(\RR^d)$ contains exactly one of them, it is called $p_{i,f}$
and its stable dimension is denoted by $n_i^s$.
\item The chart $\psi_i$ is a stable chart for $p_{i,f}$.
\end{itemize}
\end{proposition}

The major ingredient in the proof of Proposition~\ref{p=main} is the following.

\begin{prop}\label{p=diskdist}  There is a residual set $\cR_0 \subset \cC^d$ such that, for all $f\in \cR_0$,
if $x,y\in \overline{B(0,1)}\setminus \{x,y\}$ with $x\notin \cO_f(y)$, then
\begin{eqnarray}\label{e=jacdist}
\limsup_{n\to\infty} \left|\frac{\Jac(f^n)(x)}{\Jac(f^n)(y)}\right| 
= \infty.
\end{eqnarray}
\end{prop}

The proof that Proposition~\ref{p=diskdist} implies Proposition~\ref{p=main}
is quite immediate in the non-conservative case: the dynamics in any stable manifold is
diffeomorphically conjugate to a contraction of $\RR^{n^s}$; one concludes by noting that
any perturbation of the dynamics inside the stable manifold extends to a perturbation of the
dynamics on $M$. In the conservative case, this last property is much more delicate and its proof
will be postponed until Appendix~\ref{a=invmanif}.

\subsection{$K$-distortion and the Baire argument}\label{ss=baire}

In this subsection, we explain how to reduce the unbounded
distortion property (\ref{e=jacdist}) in Proposition~\ref{p=diskdist} to 
a property satisfied in finite time, which we call the $K$-distortion
property. Using a Baire argument,
we then reformulate Proposition~\ref{p=diskdist} in terms of 
this $K$-distortion property to obtain our main perturbation result 
(Theorem~\ref{t.distortion}).

\begin{defi} Let $B\subset \RR^d$ be a compact region and let $f:B\to
  B$ be an embedding. Given compact sets $\Lambda,\Delta \subset B$, 
we say that $\Lambda$ and $\Delta$ are 
{\em dynamically disjoint for $f$} if 
$f^n(\Lambda)\cap f^m(\Delta)=\emptyset$ for any $n,m\in\NN$.
\end{defi}

\begin{defi}Let $\Lambda,\Delta \subset B^d\setminus \{0\}$ be compact sets
that are dynamically disjoint for $f\in \cC^d$. 
We say that $\Lambda$ and $\Delta$ satisfy the 
{\em $K$-distortion property for $f$ at time $N$} if, 
for any $x\in \Lambda$, $y\in \Delta$ there exists $n\in \{0,\dots, N\}$ 
such that :
$$
\left|\frac{\mbox{\rm Jac } f^{n}(x)}{\mbox{\rm Jac }
  f^{n}(y)}\right|>K.
$$
More briefly, $\La, \De$ have the {\em $K$-distortion property for
  $f$}  if there exists an $N$ so that they 
satisfy the $K$-distortion property for $f$ at time $N$
\end{defi}

The properties of dynamical disjointness and $K$-distortion persist under
perturbations of both the diffeomorphism and the compact sets.
 
\begin{prop}\label{p.opendisj} Let  $\Lambda, \Delta \subset
  B^d\setminus \{0\}$ be dynamically disjoint for $f\in \cC^d$. 
\begin{enumerate}
\item  There exist neighborhoods
$\cU_f \subset \cC^d$ of $f$, $U_\Lambda$ of $\Lambda$, and $U_\Delta$ of
$\Delta$ such that all compact sets $\Lambda'\subset U_\Lambda$, $\Delta'\subset U_\Delta$ are dynamically disjoint for all $g\in \cU_f$;
\item suppose that $\Lambda$ and $\Delta$ satisfy the 
$K$-distortion property for $f$ at time $N$.  Then we can choose
$\cU_f, U_\Lambda, U_\Delta$ so that all compact sets 
$\Lambda'\subset U_\Lambda$, $\Delta'\subset U_\Delta$ 
satisfy the $K$-distortion property for $g\in \cU_f$ at time $N$.
\end{enumerate}
\end{prop}

\noindent{\bf Proof of Proposition~\ref{p.opendisj}.}
For the first item, it is enough to show that $g^n(\Lambda)\cap \Delta=\emptyset$
for $n\in \NN$. One considers $D$, a neighborhood of $0$ satisfying $f(\bar D)\subset D$
disjoint from $\Delta$. Since $\Lambda$ is compact, for large $N$ we have
$f^N(\Lambda)\subset D$. For $g$ close enough to $f$, one gets
$g^n(\Lambda)\subset D$ for each $n\geq N$, which implies the required property.

The second item is an easy continuity argument.
\eproof

We next reformulate Proposition~\ref{p=diskdist} in terms of $K$-distortion.

\begin{theo}[Main perturbation result]\label{t.distortion} Let
  $\Lambda, \Delta \subset B^d\setminus \{0\}$ be compact
sets that are dynamically disjoint for $f\in \cC^d$. 
Then for every neighborhood $\cU$ of $f$ in $\cC^d$, and for every $K>0$,
there exist $g\in \cU$ such that $\Lambda$ and $\Delta$ satisfy the 
$K$-distortion property for $g$.
\end{theo}

To prove that Theorem~\ref{t.distortion} implies Proposition~\ref{p=diskdist},
we employ a standard Baire argument.
Let $\cU_0$ be a countable basis of (relatively) 
open balls for the topology on $B^d$, and
for $U_1,U_2\in \cU$, let
$$O(U_1,U_2) = \{f\in\cD^d\,\vert\, \overline U_1\hbox{ and }\overline U_2\hbox{ are dynamically disjoint for }f\}.$$
Proposition~\ref{p.opendisj} implies that $O(U_1,U_2)$ is open in
$\cD^d$, and clearly:
\begin{eqnarray}\label{e=union}
\cD^d = \bigcup_{U_1,U_2 \in \cU} O(U_1,U_2).
\end{eqnarray}
Note that $O(U_1,U_2)$ is nonempty if and only if $\overline U_1$ and $\overline U_2$ are
disjoint, and henceforth any two such sets we discuss will be assumed
to be disjoint.

Given $U_1,U_2\in \cU$, open sets $V_1, V_2$ with
$\overline V_1\subset U_1$ and $\overline V_2\subset U_2$, 
and $K>0$, we define  $O(U_1,U_2,V_1,V_2, K)$ be the set of
all $f\in O(U_1,U_2)$ such that:
$$\overline V_1\hbox{ and }\overline V_2\hbox{ satisfy the }K\hbox{-distortion property at time }N,$$
for some $N>0$.

Theorem~\ref{t.distortion} immediately implies:
\begin{prop}\label{l=ouu1u2k} $O(U_1,U_2,V_1,V_2, K)$ is open and dense in $O(U_1,U_2)$.
\end{prop}

For each $U \in \cU$, let $\cV_{U}$ be a countable basis of open 
sets in $U$, consisting of sets whose closures are contained in $U$. Proposition~\ref{l=ouu1u2k} implies that 
$$\cR_{U_1,U_2} = \bigcap O(U_1, U_2, V_1, V_2, K)$$
is residual in $O(U_1,U_2)$, where the intersection is 
taken over all $V_1 \in \cV_{U_1}, V_2 \in \cV_{U_2}$ and $K\in\NN$.
Let $$\cR_0 = \bigcap_{U_1,U_2\in \cU} \left(\cR_{U_1,U_2} \cup \left(\cD^d\setminus O(U_1,U_2)\right)\right);$$
clearly $\cR_0$ is residual in $\cD^d$.

Suppose that $f\in \cR_0$.  Let $x,y\in B^d\setminus \{0\}$ such that $y\notin \cO(x)$
and let $K>0$ be given.
We show that there exists an  $n\in \NN$ such that
$$\left|\frac{\Jac(f^n)(x)}{\Jac(f^n)(y)}\right| > K.$$
Since $f$ is a contraction and $y\notin \cO(x)$,
there exist $U_1, U_2 \in \cU_0$ such that $x\in U_1, y\in U_2$, and
$\overline U_1,\overline U_2$ are dynamically disjoint for $f$.
This means that
$f\in O(U_1,U_2)$, and the definition of $\cR_0$ then
implies that $f\in \cR_{U_1,U_2}$.  
Let $V_1\in \cV_{U_1}, V_2\in \cV_{U_2}$ be neighborhoods
of $x$ and $y$, respectively.
Since $f\in O(U_1,U_2,V_1, V_2, K)$, we obtain that for some $n\in \NN$
$$\left|\frac{\Jac(f^n)(x)}{\Jac(f^n)(y)}\right| > K,$$
which completes the proof of Proposition~\ref{p=diskdist}.\eproof

\section{Proof of Theorem~\ref{t.distortion} (Main perturbation result)}\label{s=proof-main}

Before proving Theorem~\ref{t.distortion} we
introduce notations and concepts that will be
used in the whole section. We then isolate
the proof in some particular cases.  
We do this for two reasons: first, it will allow us to illustrate some
of the main ideas of the general case while avoiding serious technical issues,
and second, these special cases will be ingredients in the proof of the
general case.

\makeatletter
\newcommand{\subsectionruninhead}{\@startsection{subsection}{2}{0mm}
{-\baselineskip}{-0mm}{\bf\large}}
\newcommand{\subsubsectionruninhead}{\@startsection{subsubsection}{3}{0mm}
{-\baselineskip}{-0mm}{\bf\normalsize}}
\makeatother

\subsection{Preliminaries}\label{s.preliminaries}
For our purposes, a map that is linear near the origin is especially
easy to work with, because a linear map has constant Jacobian.

\subsubsectionruninhead{\em Linearization near the origin. }
A well-known feature of the $C^1$ topology is that a diffeomorphism may be $C^1$ approximated
by its derivative in a neighborhood of a fixed point:
we will say that $f\in \cC^d$ has  {\em a linear germ} if there
exists a  linear map $A:\RR^d\to \RR^d$ and a neighborhood $U$ of
$0$ in $B^d$ with $f(U)\subset U$ such that $f\vert_U =A\vert_U$.

By Proposition~\ref{p.opendisj}, if the two compact sets $\Lambda,\Delta\subset B^d\setminus \{0\}$
are dynamically disjoint for some contraction $f_0$, one can approximate $f_0$ by a contraction $f$
that has a linear germ $A=D_0f$ and such that $\Lambda, \Delta$ remain dynamically disjoint.

\subsubsectionruninhead{\em Bring $\La, \De$ into the linearized region. }
Once one considers a contraction $f$ having a linear germ, we show that
one can reduce the proof of Theorem~\ref{t.distortion}
to the case the contraction is a linear map $A$.

Let $U$ be a forward invariant set where
$f$ coincides with a linear map $A$
and choose some numbers $1>r_2>r_1>0$ satisfying
$B(0,r_2)\subset U$.
Since $f$ is a contraction, there exists an integer $m>0$ such that
$f^m(\La\cup\De) \subset B(0,r_1)$.
We would like to work with the sets $f^m(\La)$ and $f^m(\De)$ in
place of $\La$ and $\De$, and the following simple lemma allows us to do so.

\begin{lemma}\label{l.madjust}
Let $f \in \cC^d$, and let $\La, \De$ be dynamically disjoint for $f$.
For every $K>0$ and integer $m\geq 0$, there exists a neighborhood $\cV$ of $f$ and $K'>0$ such that,
for all $g\in \cV$, if $g^m(\La)$ and $g^m(\De)$ have the $K'$-distortion
property for $g$, then $\La$ and $\De$ have the $K$-distortion property for $g$.
\end{lemma}
Note that changing the Riemannian metric on $\RR^d$
only affects our choice of $K$ and the $C^1$-size of the neighborhood $\cU$ in $\cC^d$.
Hence, one can assume that $B^d$ is mapped into $B(0,1)$ by $A$ so that $A$ also is a contraction.

Let us assume that Theorem~\ref{t.distortion} has been proven for the
linear map. Then  any small perturbation
$g$ of $A$ in $\cC^d$ can be glued to $f$ inside $B(0,r_2)$ thanks to the following standard lemma:
\begin{lemma}\label{l=extend}
Given $f\in \cC^d$, numbers $1 > r_2 > r_1 > 0$, and  $\eps>0$, there exists 
$\eps'>0$ such that, for every 
embedding $g:B(0,r_2)\to \RR^d$ satisfying
$$\sup_{x\in B(0,r_2)}\|D_xg-D_x f\|_0 < \eps',$$
and $g(0)=0$, there exists a diffeomorphism $g'\in \cC^d$ such
that 
\begin{enumerate}
\item $\|g'-f\|_1 < \eps$
\item $g'=f$ on $B^d\setminus B(0,r_2)$
\item $g'=g$ on $B(0,r_1)$.
\end{enumerate}
\end{lemma}

\subsubsectionruninhead{\em Localize the perturbations. }
We introduce some terminology that will be used in the rest of the paper.
\begin{defi} Let $1>r_2>r_1 >0$. The {\em spherical shell} in $B^d$
of {\em outer radius}  $r_2$ and {\em inner radius} $r_1$ is the set:
$$S(r_1,r_2) = B(0,r_2)\setminus B(0,r_1).$$
The 
{\em modulus} $m(S(r_1,r_2))$ of the shell $S(r_1,r_2)$ is defined by:
$$m(S(r_1,r_2)) = \log\left(\frac{r_2}{r_1}\right).$$

Let $S(r_1, r_2)$  be a spherical shell.
We say that $x\in B^d$ is {\em inside} $S(r_1, r_2)$ if $\|x\| \leq r_1$ and
{\em outside} $S(r_1, r_2)$ if $\|x\| \geq r_2$.  The set of
points inside of a spherical shell 
$S$ is denoted by $I(S)$, and the set of points
outside of $S$ is denoted by $O(S)$. 

If $S(r_3, r_4)$ is another spherical shell
then we write $S(r_1, r_2) \prec S(r_3, r_4)$
if $r_2 \leq r_3$.  
We say that $x\in \RR^d$ is {\em in between} $S(r_1, r_2)$ and 
$S(r_3, r_4)$ if $x$ is
outside $S(r_1, r_2)$ and inside $S(r_3, r_4)$; that is, if
$r_2\leq \|x\| \leq r_3$. 
\end{defi}
In this terminology, the map $g'$ given by Lemma~\ref{l=extend}
coincides with $g$ inside the shell $S(r_1,r_2)$ 
and coincides with $f$ outside of $S(r_1,r_2)$.

\subsubsectionruninhead{\em Tidy perturbations. }
In order to compose several perturbations, we will often require the following property:
\begin{defi}
Let $B\subset \R^d$, let $f:B\to B$ be an embedding and let $X\subset B$.
We say that an embedding $g:B\to B$  is {\em tidy}
relative to $f$ and $X$ if, 
\begin{enumerate}
\item $g(x) = f(x)$, for all  $x\in B\setminus X$, and
\item $g^m(x) = f^m(x)$ for all $x\in B\setminus X$ and $m\geq 1$
such that $g^m(x)\in B\setminus X$.
\end{enumerate}
\end{defi}
Notice that if $X$ is forward invariant by $f$ (i.e. $f(X)\subseteq
X$), then $g$ is tidy relative to $f$ and $X$ if and only if $g=f$ on $B\setminus X$.

\subsection{The conformal case}\label{s=linconf}
We prove here Theorem~\ref{t.distortion}
in the particular case when $f=A$ is a conformal linear contraction:
we have $A=\alpha.I$ for some constant $\alpha\in (0,1)$ and some isometry $I$.
The main reason why the conformal case is simpler is that
conjugacy by a linear conformal map preserves the $C^1$-norm.
Also note that when $d=1$, this case {\em is} the general case.
We will prove the following more precise statement.
\begin{prop}\label{p.distortion}
Let $A\in \cC^d$ be a linear conformal contraction and $(\Lambda,\Delta,r)$
be a triple such that $\La,\De \subset B(0,1)$ are two
dynamically disjoint compact sets, disjoint from the ball $B(0,r)$.
Fix constants $\eps,K>0$. 

Then there exists $\beta\in (0,1)$ such that, for any $0 <s < r$,
there exists  $g\in\cC^d$ with the following properties:
\begin{enumerate}
\item $\|Dg-A\|_0 <\varepsilon$;
\item $\Lambda$ and $\Delta$ have the $K$-distortion property for $g$;
\item $g$ is tidy relative to $A$ and the spherical shell $S(\beta s,s)$.
\end{enumerate}
\end{prop}

\subsubsectionruninhead{\em The case $\Delta$ is a small ball. }
\label{sss.distortion-ball}
Once again, it is instructive to consider a simple case; we assume
first that $\De$ is a round ball $D_1$, contained in a fundamental domain 
for the action of $A$, so that $A^i(D_1)\cap D_1 = \emptyset$, for all 
$i\neq 0$.  Let $D_0\subset G$
be another round ball with the same center as $D_1$ and containing 
$\overline{D_1}$ in its interior,
chosen to be dynamically disjoint from $\La$. 

Let $\eta$ be a diffeomorphism of $B^d$ satisfying:
\begin{itemize}
\item[P1.] $\eta$ is the identity map on $B^d\setminus D_0$;
\item[P2.] the restriction of $\eta$ to the ball $D_1$ is an affine conformal
contraction whose fixed point is the center of $D_0$ and $D_1$; 
in particular, $\eta(D_1)$ is contained in the interior of $D_1$, and the 
Jacobian of $\eta$ in $D_1$ is a constant $\mu$ strictly less than $1$. 
\item[P3.] $\|D(A\circ \eta)-A\|_0<\varepsilon$,
and $\|D(A\circ \eta^{-1})-A\|_0 <\varepsilon$. 
\end{itemize}
It is easy to see that the distortion for one or more iterates
of $A\circ \eta$ between
$x\in \La$ and $y\in \De$ is equal to $\mu^{-1}$.  To get the
distortion greater than $K$, we perform a sequence of such
perturbations, each supported on a forward image $A^i(D_0)$. 

To this end, choose $m>0$ such that
$\mu^{-m} >K$.  For $i\geq 0$, the diffeomorphism $\eta_i = A^{i} \eta A^{-i}$
is supported on $A^{i}(D_0)$ and has distortion $\mu$ inside
$A^i(D_0)$. Furthermore, since conjugacy by a linear conformal map preserves the $C^1$ norm, 
the $C^1$-distance from $A\circ \eta_i$ to $A$ is the same for all $i\geq 0$,
and therefore less than $\varepsilon$. Now fix some integer $n\geq 0$ and let
$$g_{n,m}(x) = 
\begin{cases}
A\circ \eta_{n+i}(x) &\hbox{if } x\in A^{n+i}(D_0),\hbox{ for some }
i\in\{0,\ldots m-1\}, \\
A\circ \eta_{n+m+i}^{-1}(x) &\hbox{if } x\in A^{n+m+i}(D_0),\hbox{ for some }
i\in\{0,\ldots m-1\}, \\
A(x) &\hbox{otherwise. } 
\end{cases}$$
It is straightforward to check that the diffeomorphism $g=g_{n,m}$ satisfies
$$\|Dg-A\|_0 = \sup_{i\in\{0,\ldots,m-1\}} \|D(A\circ\eta_{n+i})-A\|_0 < \eps.$$
Moreover, for all $x\in  \La$, $y\in D_1 = \De$ we have
$$\left|\frac{\mbox{\rm Jac } g^{n+m}(x)}{\mbox{\rm Jac } g^{n+m}(y)}\right|>K.$$
Hence $\Lambda,\Delta$ have the $K$-distortion property for $g$.

Let us consider a point $x$ on the outside of the shell $S=S(\alpha^{n+2m}r,\alpha^nr)$
and $i>0$ such that $h^i(x)$ is inside the shell.
If $x$ does not belong to the orbit of $D_0$ for $A$, then $g_{n,m}$ coincides with $A$ on the orbit of $x$.
Let us assume now that $x$ belongs to $A^k(D_0)$ for some $k<n$:
a straightforward computation shows that
$g_{n,m}^{(n-k)+2m}(x)=A^{(n-k)+2m}(x)$ which is the first point of the orbit of $x$
that lies inside the shell $S$. Now, inside the shell $g_{n,m}=A$ and
since the inside of $S$ is forward-invariant under $A$, it follows that $g^{i}_{n,m}(x)=A^i(x)$.
Hence, $g_{n,m}$ is tidy with respect to the shell $S$.

If $\beta=\alpha^{2m+1}$ and if $n$ is the smallest integer such that $\alpha^n<s$,
the map $g=g_{n,m}$ is tidy relative to $A$ and the shell $S(\beta s, s)$.
This ends the proof of Proposition~\ref{t.distortion} in this case.

\subsubsectionruninhead{\em Cover $\Delta$ with small balls. }\label{sss.cover-ball}
For arbitrary $\Lambda, \Delta$, the strategy (to which we will return in later arguments)
is to create distortion between $\La$ and $\De$ in small increments.
Each increment will consist of a perturbation supported on a spherical
shell and will produce distortion between $\La$ and
a small piece of $\De$.

Notice that for every $r\in (0,1]$, the spherical shell
$S(\alpha r, r)$ is a fundamental domain for the action of $A$;
to simplify notations, in this section it will be denoted by $G_r$.
The following construction is an easy consequence of the fact that 
$\La$ and  $\De$ are dynamically disjoint for $A$.

\medskip
{\em There is a family $\{(D_0^j, D_1^j, r_j)\}_{j\in\{1,\dots,k\}}$ with the following properties:
\begin{itemize} 
\item for every $j$, $r_j$ is a number in $(0,1)$,
$D_0^j$ is a round disk contained in the fundamental domain $G_{r_j}$,
and $D_1^j$ is a round disk, centered at the same point as $D_0^j$ and contained in the interior of $D_0^j$;
\item $\De$ is contained in $\bigcup_{j=1}^k D_1^j$;
\item for every $j$, $\La$ is disjoint from the orbit $\bigcup_{n\in\ZZ} h^n(D_0^j)$. 
\end{itemize}}

\subsubsectionruninhead{\em Create distortion between $\Lambda$ and balls in the cover of $\Delta$. }
\label{sss.distortion-cover}
We have seen that for any ball $D_0^j$ in the cover
of $\De$, there is a perturbation $g_j$ of $A$ producing $K$-distortion
between points in $\La$ and points in $D_1^j$. These perturbations 
can be chosen to be supported in disjoint spherical shells $S_j$.
We will use the fact that $g_j$ is tidy relative to $A$ and $S_j$
in order to ensure that these perturbations can be considered independently.

\bigskip
Consider an integer $\ell>0$ such that the shell
$S(\alpha^\ell,1)$ (composed of $\ell$ successive fundamental domains of $A$)
contains $\Lambda\cup \Delta$.
Let $\beta_j$ be the constants associated to the triples $(D^j_0,D^j_1, r^j)$
by the proof of Proposition~\ref{p.distortion} given at Section~\ref{sss.distortion-ball}.
Fix now a sequence $s_1,\dots,s_k$ such that $s_1=\min(s,\alpha^{\ell})$, and 
$s_{j+1}=\alpha^{2\ell+1}\beta_j s_j$, for every $j\in\{1,\dots, k-1\}$.
Corresponding to this sequence of numbers is a sequence of shells
$S_j = S(\beta_{j} s_{j}, s_{j})$, nested as
$$S_k \prec S_{k-1} \prec \cdots 
\prec S_1 \prec  S(\alpha^\ell,1).$$
By our choice of $s_j$, in between any two successive shells  
$S_{j+1}$ and $S_j$, there
are $2\ell+1$ successive fundamental domains
$$G_{\alpha^{2\ell} \beta_j s_j}\prec\cdots\prec G_{\alpha \beta_j s_j}  \prec G_{\beta_j s_j}.$$

For each $j\in\{1,\ldots,k\}$ we denote by $g_j$ the
diffeomorphism constructed at Section~\ref{sss.distortion-ball}
for $(D_0^j, D_1^j)$ and the shell $S_j$. 
Now let $g$ be the map defined as follows:
$$g(x) = 
\begin{cases}
g_j(x) &\hbox{if } x\in S_j,\hbox{ for some }
j\in\{0,\ldots k\}, \\
A(x) &\hbox{otherwise. } 
\end{cases}$$
Then $\|Dg-A\|_0<\eps$. Moreover, one easily proves by induction on $j$ the fact that if
$x\in B^d \setminus B(0, s_1)$ and $A^n(x)$ is outside the shell $S_1$,
or between $S_{j+1}$ and $S_j$, or inside the shell $S_k$, then $g^n(x)=A^n(x)$.
In particular, $g$ is tidy with respect to $A$ and the shell $S(\beta_k s, s)$
where $\beta=\alpha^{(2\ell+1).k}\prod_{j=1}^k\beta_j$.

It remains to prove that $\La$ and $\De$ satisfy the $K$-distortion
property for $g$. Recall that $\Lambda\cup\Delta\subset S(\alpha^\ell,1)$
and that there are $2\ell +1$ fundamental domains between
any two shells $S_{j}$ and $S_{j-1}$. Hence, for every $j\in\{2,\dots, k\}$ 
there exists $n_j>0$ such that $A^{n_j}(\La\cup \De)$ 
is between the shells $S_j$ and $S_{j-1}$.
In particular $g^{n_j}=A^{n_j}$ on $\La\cup\De$. We also define $n_1=0$.

Consider $x\in \La$ and $y\in \De$ . Then there exists $j$
such that  $y\in D_1^j$. By assumption on $g_j$
there exists $N_j\in \{n_j,\dots,n_{j+1}\}$ such that
$$\left|\frac{\mbox{\rm Jac } A^{N_j}}{\mbox{\rm Jac } g_j^{N_j}(y)}\right|>K.$$
\begin{affi} For every $m\in\{n_j,\dots, N_j\}$ one has $g^{m}= A^m$
in a neighborhood of $x$.
\end{affi}
\begin{demo}
The map $g_j$ coincides with $A$ on the complement of
the orbit of $D_0^j$, which is disjoint from $\La$.
Hence, as long as $A^m(x)$ does not belong to $S_{j+1}$, one has
$g^m = A^m$ in a neighborhood of $x$. It remains to show that
$A^{N_j}(x)$ does not belong to the ball $B(0,s_{j+1})$.
As there are $2\ell+1$ fundamental domains between the shells 
$S_{j+1}$ and $S_j$, if $A^{N_j}(x)\in B(0,s_{j+1})$,
then the point $A^{N_j}(y)$ cannot belong to the shell $S_j$.
As $g_j$ is tidy with respect to $S_j$, this contradicts
the fact that  ${\mbox{\rm Jac } g_j^{N_j}(y)}\neq {\mbox{\rm Jac } A^{N_j}}$.
This contradiction concludes the proof of the claim.
\end{demo}

As a direct consequence of this claim, we obtain:
$$\left|\frac{\mbox{\rm Jac } g^{N_j}(x)}{\mbox{\rm Jac } g^{N_j}(y)}\right|>K.$$ 
Hence $\La$ and $\De$ satisfy the $K$-distortion property for $g$.

\subsection{The generic linear case, avoiding a codimension $1$ or $2$ submanifold}
\label{s.generic}
The arguments in the previous section do not generalize immediately to
the case where the linear contraction $A$ is not conformal, since
conjugation by a nonconformal linear map does not preserve the $C^1$ norm.
Nonetheless, by a small perturbation in $\cC^d$,
one can assume that the linear contraction $A$ is \emph{reduced}:
it has simple spectrum, implying that $\RR^d$ splits into a direct sum
of invariant subspaces on which $A$ is conformal.
We will denote by $F$ the $1-$ or $2-$dimensional invariant space of $A$
corresponding to the eigenvalues of $A$ of smallest modulus,
and we will denote by $E$ the sum of the other eigenspaces.

For such a reduced linear contraction, we can use an
inductive argument and the result of Section~\ref{s=linconf}.
Difficulties arise because dynamical disjointness is not preserved by projectiononto invariant subspaces. In this subsection, we treat another
special case of Theorem~\ref{t.distortion}:
the linear contraction $A$ is a reduced and the sets $\La, \De$ avoid
the codimension $1$ or $2$ submanifold $E$
determined by the weakly contracting eigenspaces of $A$.

We write $A = (A_E, A_F)$, as a product of linear contractions of $E$ and $F$:
since $A$ has simple spectrum, $A_F$ is conformal.
We denote by $U_E$ and $U_F$ the unit balls of the spaces $E$ and $F$
and by $\pi_F\colon \mathbb{R}^d \to F$ the linear projection of the product $E\times F$ on its second factor.
Note that we can change the Euclidean norm on $\RR^d$ and assume that $E$ and $F$ are orthogonal.

\subsubsectionruninhead{\em The case $\pi_F(\Lambda),\pi_F(\Delta)$ are dynamically disjoint. }
\label{sss.projection}
In this case, since $A_F$ is
conformal, Proposition~\ref{p.distortion} implies that there is a perturbation
$g_F\in \cC^{\hbox{\tiny dim}(F)}$ of $A_F$, 
tidy with respect to $A_F$ and a shell $S_F\subset U_F$, 
so that $\pi_F(\La)$ and $\pi_F(\De)$ have the
$K$-distortion property for $g_F$. Then, since
the linear map $A_E$ has no distortion, the sets
$\La$ and $\De$ also have the $K$ distortion property  
for the product map $g = (A_E, g_F)$. The embedding $g$ 
is tidy with respect to $A$ and $U_E\times S_F$ and satisfies the
conclusions of Theorem~\ref{t.distortion}.

\subsubsectionruninhead{\em Decomposition of $\Lambda$ and $\Delta$. }
In general, even though $\La$ and $\De$ are dynamically
disjoint for $A$, their projections $\pi_F(\La)$ and $\pi_F(\De)$
might not be dynamically disjoint for $A_F$.
A na\"ive way to try to fix the proof is to first
perturb $A\vert_{U_E\times U_F}$  so that the projections $\pi_F(\La)$ and $\pi_F(\De)$
become dynamically disjoint. There are two problems with this approach:
\begin{enumerate}
\item the projections $\pi_F(\La)$ and $\pi_F(\De)$ might in fact
coincide, so that a very large perturbation of $A$ would be required to disjoint them;
\item any perturbations that change the relative position of
$\pi_F(\La)$ and $\pi_F(\De)$ will destroy the invariance of
the splitting $\RR^d = E\oplus F$ in $U$, not to
mention the linear conformality of the projection $A_F$.
\end{enumerate}

To cope with these difficulties, we will cover $\La$
and $\De$ with finitely many compact, dynamically-defined pieces
$\Lambda^1,\ldots, \Lambda^{k_1}$ and
$\Delta^1,\ldots, \Delta^{k_2}$.
We perform a sequence of perturbations; at each step we
arrange for one of the pairs $(\Lambda^i, \Delta^j)$ to haven
the $K$-distortion property. Recall that in the conformal case,
we chose pieces in our cover of $\De$ to be round balls; in this
case, pieces in our covers of $\La$ and $\De$ 
will be of the form $\overline{D\times S}$, where $D$ is ball in $U_E$
and $S$ is a spherical shell in $U_F$.

\bigskip
We introduce more notation. Let $\alpha$ be the norm of the linear conformal contraction $A_F$ of $F$.
For $r\in (0,1]$, we denote by $G_r$ the shell $S(\alpha r, r)\subset U_F$,
which is a fundamental domain for the action of $A_F$ in $U_F$;
its modulus is $\mu = -\log\alpha$.

For $k\geq 1$, we define a family of spherical shells $S_{i,k}\subset U_F$,
indexed by integers $i\geq 0$, by:
$$S_{i,k} = S(\alpha^{(i+1)/k}, \alpha^{i/k}).$$
Notice that  $S_{i,k}$ has modulus $\mu/k$, for all $i, k$, and
that
$G_{1} = \bigcup_{i=0}^{k-1} S_{i,k}$.
Notice also that $A_F^j(S_{i,k})= S_{i+kj,k}$,
and so the partition
$$U_F \setminus \{0\} = \bigcup_{i = 0}^\infty S_{i,k}$$
is forward invariant under $A_F$.  

The proof of Theorem~\ref{t.distortion} in the case we are considering here can be reduced to
the following result.
\begin{prop}\label{p=4}
Let $A$ be a reduced linear contraction,
$\RR^d = E\oplus F$ its associated decomposition
and a constant $\varepsilon>0$.
Then, there exists $k_0\in \NN$ such that the following holds.

For some $u,w\in U_E$, $\delta >0$, $k\geq k_0$ and $i,j\in \NN$,
we consider the sets 
$$\La = \overline{B_E(u, \delta)\times S_{i, k}}\quad\hbox{ and }\quad 
\De = \overline{B_E(w, \delta)\times S_{j, k}},$$
and assume that
$\hat \La = \overline{B_E(u, 2\delta)\times S_{i, k}}$ and
$\hat \De = \overline{B_E(w, 2\delta)\times S_{j, k}}$
are contained in $B^d\setminus \{0\}$ and are dynamically disjoint.

Then, for any $K>0$, and any neighborhood $V_F$ of $0$ in $F$,
there exists a shell $S\subset V_F$ and an embedding
$g\colon U_E\times U_F\to U_E\times U_F$ such that:

\begin{enumerate}
\item $\|Dg - A\|_0 <\eps$;
\item $g$ is tidy with respect to $A$ and $U_E\times S$;
\item $\La,\De$ have the $K$-distortion property for $A$.
\end{enumerate}
\end{prop}

\medskip
Let us show how this proposition implies Theorem~\ref{t.distortion}.
First, Proposition~\ref{p=4} associates to $\varepsilon>0$ an integer $k_0$.
Using that $\La,\De$ are dynamically disjoint and avoid the space $E$,
one then obtains easily the following construction:

\medskip
{\em There exist $\delta>0$ and two families
$\{\La^1,\dots,\La^{k_1}\}$ and $\{\De^1,\dots,\De^{k_2}\}$
with the following properties:
\begin{itemize}
\item for every $i,j$, there exist $u_i,w_j\in E$ and $m_i,\ell_j\in \NN$ such that
$$\La^i= \overline{B_E(u_i, \delta)\times S_{   m_i, k_0}}\quad \hbox{ and } \quad
  \De^j= \overline{B_E(w_j, \delta)\times S_{\ell_j, k_0}}.$$
\item $\Lambda$ is contained in $\cup_{i=1}^{k_1}\La^i$ and $\De$ in $\cup_{j=1}^{k_2}\De^j$;
\item for every $i,j$, $\hat \La^i= \overline{B_E(u_i, 2\delta)\times S_{   m_i, k_0}}$
and $\hat \De^j=\overline{B_E(w_j, 2\delta)\times S_{\ell_j, k_0}}$
are contained in $B^d\setminus \{0\}$ and are dynamicaly disjoint.
\end{itemize}}

Fix an integer $\ell>0$ such that the shell $S(\alpha^\ell, 1)$
contains $\pi_F(\La\cup\De)$. We order all the possible pairs $(\La^i,\De^j)$
as a list $(P_1,\dots,P_{k_1k_2})$
and apply Proposition~\ref{p=4} inductively for each pair $P_m$.
We obtain a sequence of embeddings $g_m$ which are tidy with respect to $A$ and domains
$U_E\times S_m$. Since the shells $S_m$ can be chosen in arbitrarily small neighborhoods of $0$,
one can assume that they are nested as
$$S_{k_1k_2}\prec S_{(k_1k_2)-1}\prec \dots \prec S_1\prec S(\alpha^\ell,1),$$
and that between any two successive shells $S_{m+1}$ and $S_m$, there are $2\ell+1$
successive fundamental domains.

One can define the contraction $g$ and end the proof of Theorem~\ref{t.distortion}
by gluying the $g_m$ according to the domains
$U_E\times S_m$ as in Section~\ref{sss.distortion-cover}.

\subsubsectionruninhead{\em Move thin shells. }
Before proving Proposition~\ref{p=4},
we need to prove the following lemma about perturbations inside the space $F$.

\begin{lemma}\label{l=phiF} Given $\eps_0>0$, there exists $k_0\in \NN$
such that, for all  $k\geq k_0$, and for every $i\geq 2k$, 
there is a diffeomorphism $\psi:U_F\to U_F$
such that:
\begin{enumerate}
\item $\psi(\overline S_{i,k})\cap
\left(\overline S_{i+1,k}\cup \overline S_{i,k}
\cup \overline S_{i-1,k}\right)= \emptyset$;
\item $\psi = Id$ in the complement of a fundamental domain $G_r$, which is a shell
containing $S_{i+1,k}\cup \overline S_{i,k} \cup \overline S_{i-1,k}$;
\item $\|D\psi-Id_F\|_0 <\eps_0$.
\end{enumerate}
\end{lemma}

To construct $\psi$, we will use the following lemma, which will be
used again in the next subsection.

\begin{lemma}\label{l.shell}
Let $\mu_0 > 0$, $s_0\geq 2$ and $\eps_0>0$ be given.   
Then there exists $\xi>0$ such that
for any collection of conformal linear maps $f_0, f_2, \ldots f_{s_0} \in \cD^d$ satisfying
$$\sup_{i=1,\ldots, s_0}\|f_i- f_0\| < \xi,$$ 
and any collection of spherical shells $S^1 \prec S^2 \prec \cdots \prec S^{s_0} \subset B^d$
of modulus at least $\mu_0$, there exists an embedding $\psi\in \cD^d$ such that:
\begin{enumerate}
\item $\|D\psi - f_0\|_0  < \eps_0$,
\item $\psi(x) = \begin{cases}
f_0(x) &\hbox{ if } x \hbox{ is inside }S_{1},\\
f_{i}(x) &\hbox{ if } x \hbox{ is between } S_{i}\hbox{ and }S_{i+1},\hbox{ for }i=1,\ldots, s_0-1,\\
f_{s_0}(x) &\hbox{ if } x \hbox{ is outside }S_{s_0}.\\
\end{cases}$
\end{enumerate}
\end{lemma}
\begin{proof}
We prove it in the case $s_0=1$; the general case is obtained similarly.
Let $\lambda = e^{\mu_0 }>1$
and consider a smooth function $\rho\colon [0,+\infty) \to[0,1]$
which is $0$ on $[0,1]$ and $1$ on $[\lambda, +\infty)$ and whose derivative is bounded by $2/(\lambda-1)$.

Let $f_0, f_1\in \cD^d$ be two conformal linear maps.
For $r\in (0,\lambda^{-1})$, consider $\psi\in \cD^d$ defined by
$x\mapsto (1-\rho(\|x\|/r))f_0(x)+\rho(\|x\|/r) f_1(x)$; that is,
$$\psi(x)= f_0(x)+\rho\left(\frac{\|x\|}{r}\right)\left(f_1(x)-f_0(x)\right).$$
This map coincides with $f_0$ on the ball $B(0,r)$
and with $f_1$ on the complement of the ball $B(0, \lambda r)$.
Since the modulus of $S_1$ is equal to $\mu_0=\log(\lambda)$,
one can choose $r$ such that $S_1=S(r,\lambda r)$.
Moreover notice that 
$$ \|D\psi - f_0\|_0 \leq \left(1 + \frac{2\lambda}{\lambda-1}\right)\|f_1-f_0\|.$$
If one sets
$$\xi = \eps_0 \left(1 + \frac{2\lambda}{\lambda -1}\right)^{-1}$$
one thus gets $\|D\psi-f_0\|_0<\varepsilon_0$ and this completes the proof.
\end{proof}

\bigskip
We now give the proof of the first lemma.
\medskip\\
\begin{demo}[Proof of Lemma~\ref{l=phiF}]
Observe that for any $k$,
 the conformal dilation $d_k:F\to F$ defined
by $d_k(v) = \alpha^{-3/k} v$ moves the
spherical shell $S_{i,k}$ to a shell that is disjoint 
from the union $\overline S_{i+1,k}\cup \overline S_{i,k}\cup
\overline S_{i-1,k}$.  For a given $i,k$,
the map $\psi$ we construct will coincide with $d_k$
on the set $\overline{S_{i,k}}$ and with the identity outside of a fundamental
domain $G_r$ containing $S_{i,k}$ in its interior.  We will choose
 $k_0$ so that for $k\geq k_0$, the
distance from $d_k$ to the identity is small.  To insure that the
distance from $\psi$ to the identity does not depend on $i$, 
we choose the fundamental domain $G_r$ according to the following simple lemma.

\begin{affi}\label{a.bigdomain} For $k\geq 5$ and $i\geq 2k$, 
there is a fundamental domain
$G_r\subset B(0,\alpha) \subset U_F$ 
such that $\overline S_{i,k} \subset \hbox{int}(G_r)$,
and the complement of $S_{i,k}$ in $G_r$ is the union of $2$ disjoint shells:
$$G_r\setminus S_{i,k} = S^{a}\cup S^{b},$$
where $S^{a}\prec S^{b}$, and setting $\mu = -\log\alpha$ we have
$$\hbox{min}\{m(S^{a}), m(S^{b})\} \geq \mu/3.$$
\end{affi}
\begin{demo} Fix $k\geq 5$ and let
$i_0$ be an integer satisfying $4k/3 \leq i_0 \leq 5k/3-1$.
The fundamental domain $G_{\alpha} \subset U_F$ 
contains $S_{i_0,k}$, and it is easily checked
that the two shells $S^{a,k}\prec S^{b,k}$ defined by
$$G_{\alpha} \setminus S_{i_0,k} = S^{a,k}\cup S^{b,k}$$
have modulus at least $\mu/3$.
For any $i \geq 2k$ we have $i\geq i_0$ and the shell
$S_{i,k}$ is the image of $S_{i_0,k}$ under a linear conformal contraction;
the images of $G_{\alpha}, S^{a,k}$ and $S^{b,k}$ under this conformal map 
satisfy the conclusions of the lemma for $S_{i,k}$.
\end{demo}

Let $\xi>0$ be the constant specified by Lemma~\ref{l.shell},
for $\mu_0 =\mu/3$, $s_0=2$ and $\varepsilon_0$.
Choose $k_0\geq 5$ so that
$$\|d_{k_0} - Id\| = |\alpha^{-3/k_0} - 1| <\xi.$$ 
Let $i\geq 2k$, and let $G_r$ be the fundamental domain given by
the claim above. The complement of $S_{i,k}$ in
$G_r$ is the union of two spherical shells $S^a \prec S^b$ of modulus at least $\mu/3$.
Applying Lemma~\ref{l.shell} to the shells $S^a$ and $S^b$
and the maps $f_0 = f_2 = Id$ and $f_1 = d_{k_0}$, we obtain a map
$\psi$ that satisfies the conclusions of Lemma~\ref{l=phiF}.
\end{demo}

\subsubsectionruninhead{\em Separate $\Lambda$ and $\Delta$. }
We introduce a perturbation $g_1$ of $A$ that will make
some forward iterates of $\La$ and $\De$ have
dynamically disjoint $\pi_F$-projections;
this perturbation will take a special form, which will allow us to
make the final perturbation of Proposition~\ref{p=4} tidy.
\begin{lemma}\label{l.skewtype}
Let $A$ be a reduced linear contraction, $\varepsilon>0$ a constant
and $\Lambda,\Delta\subset B^d\setminus \{0\}$
two sets satisfying the assumptions of Proposition~\ref{p=4}.
 
Then, for any neighborhood $V_F$ of $0$ in $F$,
there exists a diffeomorphism $g_1=A\circ \varphi$ of $\RR^d$, a shell $S_1\subset V_F$ of $F$
which is a fundamental domain of $A_F$, and an integer $\ell_1\geq 1$ such that:
\begin{enumerate} 
\item the map $\varphi$ takes the form
$\varphi (u,v) = (u, \theta(u,v))$, for $(u,v)\in E\times F$, and
is supported in $U_E\times S_1$;
\item $\max\left(\|D\varphi-\mbox{Id}\|_0, \; \|D\varphi^{-1}-\mbox{Id}\|_0\right)
<\frac{\varepsilon}{\|A\|}$,
and thus $\|Dg_1 - A\|_0 < \eps$;
\item $g_1^{\ell_1}(\hat\La\cup\hat\De)$ is contained
in $U_E\times I(S_1)$, and 
$\pi_F(g_1^{\ell_1}(\hat\La))$ and $\pi_F(g_1^{\ell_1}(\hat\De))$ 
are dynamically disjoint for the restriction of $A_F$ to $I(S_1)$.
\end{enumerate}
\end{lemma}
\begin{proof}
Let $\varepsilon_0=\frac \varepsilon {4\|A\|}$; Lemma~\ref{l=phiF}
associates to this constant an integer $k_0$. For some integer $k\geq k_0$,
we consider two compact sets
$$\La = \overline{B_E(u, \delta)\times S_{i, k}}\quad\hbox{ and }\quad 
\De = \overline{B_E(w, \delta)\times S_{j, k}}.$$
The same formula with balls in $E$ of radius $2\delta$ defines the compact sets $\hat \Lambda,\hat \Delta$.
Recall that these two sets are assumed to be contained in
$B^d\setminus \{0\}$ and to be dynamically disjoint.

Let $\nu = \|A_E^{-1}\|\|A_F\|$; because $F$ is the maximally contracted eigenspace
of $A$, we have $\nu <1$. Fix a large integer $\ell$ such that
$$i+k\ell\geq 2k \quad \text{ and } \quad \nu^{\ell} < \frac {\delta\varepsilon_0} 4.$$
If $\pi_F(A^{\ell}(\La))$ and $\pi_F(A^{\ell}(\De))$ are
dynamically disjoint for $A_F$, then there
is nothing to prove (we choose $S_1$ to be any spherical shell, set $g_1=A$
and choose $\ell_1\geq \ell$ such that $g_1^{\ell_1}(\hat\La\cup\hat\De)$ is contained
in $U_E\times I(S_1)$).

If, on the other hand, $\pi_F(A^{\ell}(\La))=\overline{S_{i+k\ell,k}}$
and $\pi_F(A^{\ell}(\De))=\overline{S_{j+k\ell,k}}$ are not
dynamically disjoint for $A_F$, then there exists $q_0\in \ZZ$ such that
$$A_F^{\ell}\left(\overline{S_{i,k}}\right)\cap A^{q_0}_F\left(\overline{S_{j,k}}\right)\neq\emptyset.$$
This implies that for some $\kappa\in\{-1,0,1\}$ we have
$$A^{q_0}_F(\overline{S_{j,k}}) = \overline{S_{i+k\ell+\kappa,k}}.$$

Since $i+k\ell\geq 2k$, Lemma~\ref{l=phiF} for $\varepsilon_0$ and the shell $S_{i + k\ell, k}$
provides us with a fundamental domain $S_1=G_r$
and a diffeomorphism $\psi\in \cD^{{\hbox{\tiny dim}}(F)}$ such that
\begin{enumerate}
\item $S_1 \cap  \bigcup_{q\in \NN} A^q_F(\overline{S_{j, k}})=A^{q_0}_F(\overline{S_{j, k}})
=\overline{S_{i+k\ell+\kappa,k}}$,
\item $\psi(\overline S_{i+k\ell,k})\cap \left(\overline S_{i + k\ell +1,
k}\cup \overline S_{i+k\ell, k}\cup \overline S_{i + k\ell-1,k}\right)= \emptyset$,
\item $\psi = Id_F$ in the complement of $S_1$, and
\item $\|D\psi-Id_F\|_0 < \eps_0$.
\end{enumerate}

Fix a smooth bump function $\rho_0:U_E\to [0,1]$ which is $1$ on $B(u,\delta)$
and $0$ outside $B(u,2\delta)$, and whose derivative has a norm bounded by $2\delta^{-1}$.
One then defines $\rho: A_E^\ell(U_E)\to [0,1]$ and a map $\theta : A^{\ell}(\hat\La)\to S_1$ by:
\begin{eqnarray*}
\rho(u) &=& \rho_0(A_E^{-\ell}(u)),\\
\theta(u,v) &=& \rho(u)\psi(v) + (1-\rho(u)) v.
\end{eqnarray*}
Finally, for $(u,v)\in A^{\ell}(\hat\La)$, let
$$\varphi(u,v)= (u,\theta(u,v));$$
we extend this definition by the identity to obtain a diffeomorphism $\varphi\in \cD^d$.
It is supported in $A^{\ell}(\hat\La)\subset U_E\times S_1$ and coincides with $\mbox{Id}\times \psi$
on $A^\ell(\La)$.

Note that the support of $\psi$ is contained in the ball $A_F^\ell(U_F)$ with radius $\|A_F\|^\ell$.
Hence, we have
\begin{eqnarray*}
\|D\varphi - Id\|_0 &\leq& \|D\rho\|_0\,\|\psi -
Id_F\|_0 + \|\rho\|_0\,\|D\psi - Id_F\|_0\\
&<& 2\delta^{-1}\,\|A_E^{-1}\|^\ell\cdot2\|A_F\|^\ell + \eps_0\\
&<& 4 \delta^{-1} \nu^\ell +  \eps_0\\
&<& 2\eps_0=\frac{\varepsilon}{2\|A\|}.
\end{eqnarray*}
An elementary calculation shows that, for any diffeomorphism $h$ of $B^d$,
if $\|Dh - Id\|_0 \leq \frac{1}{2}$, then $\|Dh^{-1} - Id\|_0 \leq 2 \|Dh - Id\|_0$.
Hence, we get
\begin{eqnarray}\label{e=varphinv}
\|D\varphi^{-1} - Id\|_0 &\leq& 2\|D\varphi - Id\|_0\, \leq\, \frac\eps{\|A\|}.
\end{eqnarray}

Let $g_1 = A\circ \varphi$ and choose $\ell_1 > \ell$
so that $g_1^{\ell_1}(\hat\La\cup\hat\De)$ is contained in $U_E\times I(S_1)$.
Since $\varphi$ is supported in $A^{\ell}(\hat\La)$, and since $\hat \La$
and $\hat \De$ are dynamically disjoint for $A$, we have for each $k\geq 0$,
$g_1^{k}(\De)=A^{k}(\De)$.
Since
$$S_1 \cap  \bigcup_{q\in \ZZ} \pi_F \left(A^q(\hat\De)\right)
= \pi_F(A^{q_0}(\hat\De)) \subset \overline S_{i + k\ell +1,k}\cup
\overline S_{i+k\ell, k}\cup 
\overline S_{i + k\ell-1,k},$$
and
$$\pi_F(A^{\ell}(\hat\La))=\overline S_{i+k\ell, k},$$
it follows from the properties of $\psi$ that 
$\psi(\pi_F(A^{\ell}(\La)))$ and $\pi_F(A^{q_0}(\De))$ are
dynamically disjoint for $A_F$.
Since $g_1=A$ on $U_E\times S_1$ and $g_1=A\circ \varphi$
on $A^\ell(\hat \La)$, this immediately implies that 
$\pi_F(g_1^{\ell_1}(\La)$ and $\pi_F(g_1^{\ell_1}(\De))$ are
dynamically disjoint for the restriction of $A_F$ to $I(S_1)$.
This completes the proof of Lemma~\ref{l.skewtype}.
\end{proof}

\subsubsectionruninhead{\em Proof of Proposition~\ref{p=4}. }
The desired perturbation $g$ of $A$ will be obtained in three steps $g_1,g_2,g_3$.
We first choose an integer $k_0\in \NN$ according to Lemma~\ref{l=phiF}.
We apply Lemma~\ref{l.skewtype} and obtain the first perturbation $g_1$
supported on a set $U_E\times S_1$ where $S_1$ is a shell in $F$ contained in an arbitrarily
small neighborhood of $0$. This also provides us with two dynamically disjoint
sets $\pi_F(g_1^{\ell_1}(\Lambda))$ and $\pi_F(g_1^{\ell_1}(\Delta))$ for the dynamics of $A_F$.

We then apply the argument of Section~\ref{sss.projection}: fixing a constant $K'>0$,
we obtain a perturbation $h=(A_E,g_F)$ of $A$ and a shell $S_2\prec S_1$ of $F$
such that $g_1^{\ell_1}(\Lambda)$ and $g_1^{\ell_1}(\Delta)$ have the $K'$-distortion property
for $h$. Moreover, $h$ is tidy with respect to $A$ and $U_E\times S_2$.
We define $g_2$ as the map which coincides with $g_1$ on $B^1\setminus (U_E\times S_2)$
and with $h$ on $U_E\times S_2$. 
If $K'$ has been chosen large enough, then $(\Lambda,\Delta)$ will 
have the $K$-distortion property at some time $N>\ell_1$ for $g_2$.

In our final perturbation, we go from the untidy map $g_2$
to a map $g_3$ that is tidy with respect to
a larger spherical shell containing $S_1$ and $S_2$, while at the same
time keeping the desired distortion properties of $g_2$.
We choose $\ell_2> N$ so that 
$A_F^{\ell_2}(\bar S_1)\subset \hbox{int }I(S_2)$.
Recall that there exists a diffeomorphism $\varphi$ supported on $U_E\times S_1$ such that
$g_1=A\circ \varphi$. We define $g_3=g_2\circ \bar \varphi$ where,
$$\bar \varphi=A^{\ell_2}\circ \varphi^{-1}\circ A^{-\ell_2}.$$
If one sets $S_3=A_F^{\ell_2}(S_1)$, one sees that
the diffeomorphisms $g_3$ and $g_2$ coincide outside the set $U_E\times S_3$, which is contained in a small
neighborhood of $E$, if $\ell_2$ has been chosen large enough.
In particular, the sets $(\Lambda,\Delta)$ have the $K$-distortion property for $g_3$.

The maps $g_3$ and $A$ coincide outside three domains.
On $U_E\times S_1$ and $U_E\times S_2$, we have $\|Dg_3-A\|_0<\varepsilon$.
On $U_E\times S_3$, the claim below gives
$$\|Dg_3-A\|_0<\|A\|.\|D\bar \varphi-\mbox{Id}\|_0<\|A\|.\|D\varphi^{-1}-\mbox{Id}\|_0,$$
which is less than $\varepsilon$ by Lemma~\ref{l.skewtype}.

\begin{affi} 
We have $\|D\bar \varphi-\mbox{Id}\|_0<\|D\varphi^{-1}-\mbox{Id}\|_0$.
\end{affi}
\begin{demo}
Since $\varphi^{-1}$ has the form $\varphi^{-1}(u,v)=(u,\bar \theta(u,v))$, we have,
$$ D \varphi^{-1} = 
\left(
\begin{array}{cc}
Id_E & 0 \\\overline\theta_u  &
\overline\theta_v
\end{array}
\right).$$
This gives
\begin{eqnarray}\label{e=dphiinv}
D \varphi^{-1} = 
\left(
\begin{array}{cc}
Id_E & 0 \\A_F^{\ell_2}\overline\theta_uA_E^{-\ell_2}  &
A_F^{\ell_2}\overline\theta_v A_F^{-\ell_2}
\end{array}
\right).
\end{eqnarray}
Since $\|A_F\|\|A_E^{-1}\|<1$, we have
$\|A_F^{\ell_2}\overline\theta_uA_E^{-\ell_2}\|\leq \|\overline\theta_u\|$.
Using that $A_F$ is conformal and since $E,F$ were assumed to be orthogonal,
one gets the announced inequality.
\end{demo}

Let $S$ be the smallest shell in $F$ that contains $S_1$, $S_2$ and $S_3$.
Note that all these shells can be constructed in any small neighborhood $V_F$ of $0$ in $F$.
It remains to prove that $g=g_3$ is tidy with respect to $A$ and $U_E\times S$.
Consider a point $x\in U_E\times U_F$ and an integer $m$ such that $\pi_F(x)$
is outside $S$ and $\pi_F(g_3^m(x))$ is inside $S$. Since $S_1$ and $S_3$
are fundamental domains of $A_F$, there exist unique integers $i_1$ and $i_3$
such that $g_3^{i_1}(x)\in U_E\times S_1$ and $g_3^{i_3}(x)\in U_E\times S_3$;
moreover $i_3=i_1+\ell_2$. Thus,
$$g^m_3(x)=A^{m-i_1-\ell_2}\bar \varphi h^{\ell_2}\varphi A^{i_1}(x).$$
Since $h$ is tidy with respect to $A$ and $U_E\times S_2$,
and by definition of $\bar \varphi$, we have
$$g^m_3(x)=A^{m-i_1-\ell_2}A^{\ell_2}\varphi^{-1} A^{-\ell_2}A^{\ell_2}\varphi A^{i_1}(x)=A^m(x).$$
This ends the proof of Proposition~\ref{p=4}.

\subsection{The general case}
We now prove Theorem~\ref{t.distortion}; as we saw in Sections~\ref{s.preliminaries} and~\ref{s.generic},
we may assume without loss of generality that $f$ is a reduced linear contraction $A$.
The proof will be by induction on the dimension $d$.
Note that the case $d=1$ is a direct consequence of Proposition~\ref{p.distortion},

Let us suppose that Theorem~\ref{t.distortion} has been proved
in any dimension $d'<d$ and let $A$ be a reduced linear contraction of $\RR^d$.
As in Section~\ref{s.generic}, one introduces the associated decomposition $\RR^d=E\oplus F$.
Let $\La, \De \subset B^d\setminus\{0\}$ be two compact sets
that are dynamically disjoint for $A$, and fix a constant $K>0$.
The desired perturbation will be obtained in three steps $g_1,g_2,g_3$:
using the induction hypothesis, we will first obtain distortion between
$\La\cap E$ and $\De\cap E$; next, we will create the distortion property
between $\La\cap E$ and $\De \setminus W$ and between $\La\setminus W$
and $\De\cap E$, where $W$ is a small neighborhood of $E$ in $\RR^d$;
in the last step, we will use the results of Section~\ref{s.generic} to 
complete the proof of the theorem, obtaining the distortion between $\De\setminus W$ and
$\La\setminus W$.

\subsubsectionruninhead{\em Distortion property on the weak-stable space:
use of the induction hypothesis. }
We first consider the induced dynamics $A_E$ of $A$ on $E$:
the compact sets $\De_E= \De\cap U_E$ and $\La_E= \La \cap U_E$
are dynamically disjoint for $f_E$.

By the induction hypothesis, there exists a contraction $h_E$ of $E$, arbitrarily close to $A_E$ such that
$\De_E$ and $\La_E$ have the $K$-distortion property for $h_E$. Note that one can again perturb $h_E$
in a small neighborhood of $0$ and assume furthermore that $h_E$ coincides with $A_E$ near $0$.

We define $g_1=(h_E, A_F)$; by continuity,
there exist  compact neighborhoods
$\hat \De_E$ and $\hat \La_E$ of $\De_E$ and $\La_E$, respectively,
that satisfy the $K$-distortion property
for any embedding close to $g_1$.

\subsubsectionruninhead{\em Distortion between points in the weak stable space
and points not in the weak stable space. }
Let $\La' = \La \setminus \hbox{int }(\hat \La_E)$ and
$\De' = \De \setminus \hbox{int }(\hat \De_E)$.
We will obtain the $K$-distortion property between $\La'\cup \De'$ and $\De_E\cup\La_E$.
Since these sets belongs to disjoint forward invariant regions,
the construction will be much easier than in Section~\ref{s.generic}.

Fix a neighborhood $V_1$ of $0$ that is forward invariant and where $g_1$
coincides with $A$. We choose an integer $\ell_1$ such that $g_1^{\ell_1}(\La'\cup \De')$
is contained in $V_1$.
Working in $F$, we note that the set $\Gamma=\pi_F(g_1^{\ell_1}(\La'\cup \De'))$ is disjoint from a ball
$B(0,r_0)$. We consider any triple of positive numbers $(r_1,r_2,r_3)$ such that
the shells $S_{r_i}=S(\frac 1 2 r_i,r_i)$ are contained in $B(0,r_0)$
and satisfy $S_{r_3}\prec S_{r_2} \prec S_{r_1}$. Their modulus is equal to $\mu_0=\log 2$.
By Lemma~\ref{l.shell}, we associate to $\mu_0$, $s_0=4$ and a small $\varepsilon_0$, a constant $\xi>0$.

Choose a number $\beta > 1$ and conformal dilations $f_1, f_2$ of $F$ such that
$$\|f_1 - A_F\| < \xi, \;\; \|f_2 - A_F\| < \xi, \;\;
\left|\frac{\det f_1}{\det A_F}\right| > \beta,\;\;
\hbox{ and } \left|\frac{\det A_F}{\det f_2}\right| > \beta.$$
Let $\psi$ be the map associated by Lemma~\ref{l=phiF} to the $S_{r_i}$, and the maps $f_0=f_3= A_F$, $f_1, f_2$.
Then $\psi$ is a contraction which is close to $A_F$ if the constant $\varepsilon_0$ has been chosen
small enough, that coincides with $A_F$ outside $B(0, r_1)$ and in a neighborhood of $0$,
with $f_1$ between the shells $S_{r_1}$ and $S_{r_2}$, and with $f_2$ between the shells
$S_{r_2}$ and $S_{r_3}$. The constants $r_1,r_2,r_3$ are now chosen
according to the following property.

\begin{affi}
For any constant $K_0>0$, and given any $r_1$, there exists $C_1>0$ such that,
for any $r_2$ in $(0, r_1/C_1)$, and for any $r_3$, the map $\psi$ satisfies
$$\forall z\in\Gamma, \; \exists n>0,\, \;\;
\left|\frac{\mbox{\rm Jac } \psi^{n}(z)}{\mbox{\rm Jac } \psi^{n}(0)}\right|>K_0.$$

Furthermore, given any $r_1,r_2$, there exists $C_2>0$ such that,
for any $r_3$ in $(0,r_2/C_2)$, the map $\psi$ satisfies 
$$\forall z\in\Gamma, \; \exists m>0,\, \;\;
\left|\frac{\mbox{\rm Jac } \psi^{m}(0)}{\mbox{\rm Jac }\psi^{m}(z)}\right|>K_0.$$
\end{affi}
\begin{demo} When $\frac{r_1}{r_2}$ goes to infinity, the orbits of 
the points $z\in \Gamma$ spend an interval of  times 
$\{n_1(z),\dots,n_2(z)\}$ between the shells $S_{r_1}$ and 
$S_{r_2}$ whose length $n_1(z)-n_2(z)$ goes to infinity; 
furthermore, $n_1(z)$ does not depend of $r_2, r_3$ and is 
uniformly bounded on $\Ga$ by some integer $n_1$. At each 
iteration between $S_{r_1}$ and $S_{r_2}$ the distortion 
$\left|\frac{\mbox{\rm Jac } \psi^{n}(z)}{\mbox{\rm Jac } \psi^{n}(0)}\right|=\left|\frac{\mbox{\rm Jac } \psi^{n}(z)}{\mbox{\rm Jac } A_F^{n}}\right|$ increases by the 
factor $\beta$.  Hence it is enough to choose $r_2$ such that $n=\inf_{z\in\La} n_2(z)-n_1(z)$
satisfies
$$\beta^n \left|\frac{\mbox{\rm Jac }\psi^{n_1}(z)}{\mbox{\rm Jac } A_F^{n_1}}\right|>K_0.$$
The proof of the second part of the claim is analogous. 
\end{demo}

We now define $g_2=(h_E,\psi)$. 
Then, for any $x\in g_2^{\ell_1}(\La'\cup \De')$
and any $y\in g_2^{\ell_1}(\De_E\cup \La_E)$, there exists an integer $n>0$ such that
$$\frac{\mbox{\rm Jac }g_2^{n}(x)}{\mbox{\rm Jac }g_2^{n}(y)} =
\frac{\mbox{\rm Jac }\psi^{n}(\pi_F (x))}{\mbox{\rm Jac }A_F^{n}(0)}>K_0,$$
and similarly there exists an integer $m>0$ such that
$\frac{\mbox{\rm Jac }g_2^{m}(y)}{\mbox{\rm Jac }g_2^{m}(x)}>K_0$.
If $K_0$ has been chosen large enough (with respect to $\ell_1$), one deduces that
$\La'\cup \De'$ and $\De_E\cup\La_E$ have the $K$-distortion property for $g_2$.

\subsubsectionruninhead{\em Distortion between points of $\La, \De$ in the complement of the weak stable space. }
From the two previous steps, we have shown that the pairs $(\La_E, \De)$ and $(\La, \De_E)$ 
satisfy the $K$-distortion property for any embedding close to $g_2$.
Hence there are compact neighborhoods $O_\La$ and $O_\De$
of $\La_E$ and $\De_E$ respectively, such that the
pairs $(O_\La, \De)$ and $(\La, O_\De)$  satisfy the
$K$-distortion property for $g_2$.
Let $\La'' = \La\setminus \hbox{int }(O_\La)$,
and let $\De'' = \De\setminus \hbox{int }(O_\De)$.
All that remains is
to create distortion between $\La''$ and $\De''$ 
for an embedding $g_3$ close to $g_2$. 

By construction, there exists a neighborhood $V_2$ of $0$ that is forward invariant and where $g_2$
coincides with $A$. We choose some integer $\ell_2$ such that $g_2^{\ell_2}(\La''\cup \De'')$
is contained in $V_2$. It remains to apply Theorem~\ref{t.distortion} as it was proved at
Section~\ref{s.generic}, to the map $A$, the sets $g_2^{\ell_2}(\La''),g_2^{\ell_2}(\De'')$
and a constant $K_0$ large. We obtain an embedding $g=g_3$
such that $\La,\De$ have the $K$-distortion property.
The set where $g$ and $A$ differ has been constructed in an arbitrarily
small neighborhood of the origin.

\section{Large centralizer for a locally $C^1$-dense set of diffeomorphisms}\label{s=contre-exemple}
\subsection{The case of the circle}
Our aim in this subsection is to prove the first part of Theorem~\ref{t=main2}.
The following lemma summarizes some very classical properties of diffeomorphisms of the circle. 

\begin{lemm}\label{l.reduction}
Let $\cD^1_0\subset \Diff^{~1}(S^1)$ be the set of $f\in\cD^1_0$ satisfying the following properties:
\begin{itemize}
\item $f$ is a $C^\infty$ Morse-Smale diffeomorphism (i.e. the non wandering set consists of
finitely many hyperbolic periodic points, alternately attracting or
repelling); and
\item for every periodic point $x\in Per(f)$, there is a neighborhood $U_x$ of $x$ such that  the restriction  $f|_{U_x}\colon U_x\to f(U_x)$ 
is  an affine map (for the natural affine structure on $S^1=\RR/\ZZ$).
\end{itemize} 
Then $\cD^1_0$ is dense in $\Diff^1(S^1)$.   
\end{lemm}

For $\alpha>1$ and $\beta\in (0,1)$,
one introduces the set $D_{\alpha,\beta}$ of orientation preserving $C^\infty$ diffeomorphism
of the interval $[0,1]$ with the following properties:
\begin{itemize} 
\item $\{0,1\}$ is the set of fixed points of $f$, and $f(x)>x$ for $x\in(0,1)$;
\item $f(x)=\alpha x$ for small $x$ and $f(x)=1+\beta(x-1)$ for $x$ close to $1$. 
\end{itemize}

Now item 1 of  Theorem~\ref{t=main2} is a consequence of the following proposition
(see Section~\ref{sss=circle}).

\begin{prop}\label{p.contre-exemple} Let $f$ be a diffeomorphism in $D_{\alpha,\beta}$.
Then any $C^1$ neighborhood $\cU$ of $f$ in $\Diff^1([0,1])$ contains a diffeomorphism $g$
such that $g=f$ in a neighborhood of $\{0,1\}$
and $g$ is the time one map of a $C^\infty$-vector field on $[0,1]$. 
\end{prop}

\subsubsection{Mather invariant}
We recall here a construction introduced by J. Mather~\cite{mather} which associates to any diffeomorphism
$f\in D_{\alpha,\beta}$ a class of diffeomorphism of $S^1$.

Let us fix $\alpha>1$ and $\beta\in (0,1)$ and introduce a $C^\infty$ orientation preserving diffeomorphism 
$\varphi\colon (0,1)\to \R$ such that $\varphi(x)= \frac{\ln x}{\ln \alpha}$ for $x$ small and 
$\varphi(x)= \frac{\ln (1-x)}{\ln \beta}$  for $x$ close to $1$: there exists $K_0>0$ such that
$\varphi^{-1}(x)= e^{\ln\alpha. x}$ for $x<-K_0$ and  $\varphi^{-1}(x)= 1-e^{\ln\beta. x}$ for $x>K_0$.

For any $f\in D_{\alpha,\beta}$  the conjugated diffeomorphism $\theta_f=\varphi\circ f\circ \varphi^{-1}$
of $\RR$ satisfies $\theta_f(x)>x$ for all $x$; furthermore $\theta_f(x)$ coincides with $x+1$
when $|x|$ is larger than a constant $K_f>K_0$.

The space $\RR/\theta_f$ of the orbits of $\theta_f$ is a smooth circle $S_f$ which has
two natural identifications with the (affine) circle $S^1=\RR/\ZZ$:
two points $x,y\in (-\infty,-K_f]$ (resp. $x,y\in [K_f,+\infty)$) are in the  same orbit for $\theta_f$
if and only if they differ by an integer. This leads to two diffeomorphisms
$\pi_+\colon S_f\to S^1$ and $\pi_-\colon S_f\to S^1$, respectively.
Let  $\De_{f,\varphi}=\pi_+\circ \pi_-^{-1}\colon S^1\to S^1$.

\begin{lemm} \label{l.matherdef}
The diffeomorphism $f$ is the time-one map of a $C^1$ vector field if and only if $\De_{f,\varphi}$ is a rotation.
\end{lemm}
\begin{demo} Note that $f\in D_{\alpha,\beta}$ coincides with the time one map of the vector field
$X^-=\ln\alpha. x\frac\partial{\partial x}$ in a neighborhood of $0$ and 
with $X^+= \ln\beta.(x-1)\frac\partial{\partial x}$ in a neighborhood of $1$.
Furthermore, if $f$ is the time one map of a $C^1$ vector field $X$ on $[0,1]$ then $X=X^-$ in a neighborhood of $0$ and $X=X^+$ in a neighborhood of $1$. 
The hypothesis on $\varphi$ implies $\varphi_*(X_-)=\frac\partial{\partial x}$ on some interval $(-\infty,L_-)$
and $\varphi_*(X_+)=\frac\partial{\partial x}$ on an interval $(L^+,\infty)$.

Assume that $\De_{f,\varphi}$ is a rotation. Then we define a vector
field  $Y$ on $\RR$ as follows: consider $n>0$ such that
$\theta_f^n(x)> K_f$. Now let $Y(x)=\left(D_x\theta_f^n \right)^{-1}(\frac\partial{\partial x})$. This vector does not depend of $n$ (because $\theta_f$ is the translation $t\mapsto t+1$ for $t\geq K_f$). 

\begin{clai} if $x<-K_f$ then $ Y(x)=\frac\partial{\partial x}$
\end{clai}
\begin{demo} Consider the natural projection $\pi_f\colon \RR\to S_f$ that maps each point to its orbit for $\theta_f$.
Since $Y$ is invariant by $\theta_f$, the vector field $(\pi_f)_*(Y)$ is well-defined. Since on $(K_f,+\infty)$
the vector $Y(x)$ is equal to $\frac\partial{\partial x}$, the map $\pi_+\circ \pi_f$
coincides with the natural projection $\RR\to S^1$ and we have
$(\pi_+\circ\pi_f)_*(Y(x))=\frac\partial{\partial x}$.
As $\De_{f,\varphi}$ is a rotation, and as the rotations preserve the vector field $\frac\partial{\partial x}$,
we obtain that 
$(\pi_-\circ\pi_f)_*(Y(x))=(\De_{f,\varphi}^{-1}\circ \pi_+\circ\pi_f)_*(Y(x))=\frac\partial{\partial x}$.
As $\theta_f$ coincides with the translation $t\mapsto t+1$ on $(-\infty,-K_f]$, the projection $\pi_-\circ \pi_f$
coincides on $(-\infty,-K_f]$ with the natural projection $\RR\mapsto\RR/\ZZ$.
Hence $(\pi_-\circ\pi_f)_*(Y(x))=\frac\partial{\partial x}$ implies $Y(x)=\frac\partial{\partial x}$.
\end{demo} 

Notice that, by construction, the vector field $Y$ is invariant by $\theta_f$;
furthermore $\theta_f$ is the time one map of $Y$: this is true on  a neighborhood of $\pm\infty$,
and extends on $\RR$ because $Y$ is $\theta_f$-invariant. 

Now,  the vector field $X=\varphi^{-1}_*(Y)$, defined on $(0,1)$, coincides with $X_-$ and $X_+$ in a neighborhood
of $0$ and $1$, respectively, hence induces a smooth vector field on $[0,1]$.
Finally, $f$ is the time one map of $X$.

Conversely, if $f$ is the time one map of a $C^1$-vector field $X$ on $[0,1]$ then $\theta_f$ is the time one map
of  the vector field $Y=\varphi_*(X)$, which coincides with $\partial/\partial x$ in the neighborhood
of $\pm\infty$ (because $X$ coincides with $X_-$ and $X_+$ in a neighborhood of $0$ and $1$, respectively).
Hence the projections $(\pi_-\circ\pi_f)_*(Y)$ and $(\pi_+\circ\pi_f )_*(Y)$ are both equal
to the vector field $\partial/\partial x$ on $S^1$. This implies that
$(\De_{f,\varphi})_* (\partial/\partial x)=\partial/\partial x$, which
implies that $\De_{f,\varphi}$ is a rotation. 
\end{demo}

\begin{rema}
The function $\De_{f,\varphi}$ defined here seems to depend on the choice of $\varphi$.
There is a more intrinsic way to define the diffeomorphism  $\De_{f,\varphi}$ ``up to composition by rotation": 

The vector fields $X_-$ and $X_+$ defined in a neighborhood of $0$ and $1$, respectively,
are the unique vector fields such that  $f$ is the time one map of the corresponding flows,
in the neighborhood of $0$ and $1$, respectively. Each of these vector fields induces
a parametrization of the orbit space $(0,1)/f=S_f$, that is, up to the choice of an origin,
a diffeomorphism $\pi_f^\pm\colon S_f\to S^1$. The change of parametrization $\pi_f^+\circ (\pi_f^-)^{-1}$
is well defined, up to the choice of an origin of the circle, i.e. up to multiplication,
at the right and at the left, by rotations. This class of maps is
called the {\em Mather invariant} of $f$.  
\end{rema}

\subsubsection{Vanishing of the Mather invariant: proof of Proposition~\ref{p.contre-exemple}\label{ss.contre-exemple}}
Fix $f\in D_{\alpha,\beta}$ and $K_f>0$ such that $\theta_f=\varphi \circ f\circ \varphi^{-1}$
coincides with $x\mapsto x+1$ on $(-\infty,-K_f]\cup[K_f,+\infty)$. 

Given a diffeomorphism $h\colon \RR\to\RR$, the \emph{support of $h$}, denoted by $supp(h)$
is the closure  of the set of points $x$ such that $h(x)\neq x$.
\begin{lemm} \label{l.composition} Consider a number $a>K_f$ and a diffeomorphism
$\tilde \psi\colon\RR\to\RR$  whose support is contained in $(a,a+1)$.
Let $h$ denote the  diffeomorphism $\varphi^{-1}\circ \tilde
\psi\circ\varphi$, and let
$\psi$ denote the diffeomorphism of $S^1\simeq [a,a+1]/a\sim a+1$ induced by $\tilde \psi$. 

Then the diffeomorphism $g= f\circ h$ belongs to $D_{\alpha,\beta}$, and
$\De_{g,\varphi}=\psi \circ\De_{f,\varphi}$. 
\end{lemm}
\begin{demo} The diffeomorphism $g$ coincides with $f$ in neighborhoods of $0$ and $1$ proving that
$g\in D_{\alpha,\beta}$. Furthermore, by construction, one may choose $K_g=a+1$. 

If $x<-a$, there is a (unique) integer such that $\theta_f^n(x)=\theta_g^n(x)\in[a,a+1)$,
and by construction of $\De_{f,\varphi}$,  the projection of $\theta_f^n(x)$ on $S^1$ is $\De_{f,\varphi(x)}$.
Now the projection on $S^1$ of $\theta_g^{n+1}(x)=\theta_f\circ \tilde \psi\circ\theta_f^n(x)$
is $\psi\circ\De_{f,\varphi}(x)$, by construction. 
As $\theta_g=\theta_f=y\mapsto y+1$ for $y\geq a+1$, one gets that the projection on $S^1$ of $\theta_g^{n+k}(x)$
is $\psi\circ\De_{f,\varphi}(x)$, for all $k>0$; hence $\De_{g,\varphi}=\psi\circ\De_{f,\varphi}$.
\end{demo}

Iterating the process described in Lemma~\ref{l.composition} one gets:

\begin{coro}\label{c.composition} Consider a finite sequence of numbers $a_i>K_f$, $i\in\{1,\dots,\ell\}$,
such that $a_{i+1}>a_i+1$ for all $i\in\{1,\dots,\ell-1\}$. For every $i\in\{1,\dots,\ell\}$, fix
a diffeomorphism $\tilde \psi_i\colon\RR\to\RR$ whose support is contained in $(a_i,a_i+1)$.
Let $h_i$ denote the diffeomorphism $\varphi^{-1}\circ \tilde \psi_i\circ\varphi$,
and let $\psi_i$ denote the diffeomorphism of $S^1$ induced by $\tilde \psi_i$.
(Note that the diffeomorphisms $h_i$ have disjoint support, so that they are pairwise commuting.)

Then the diffeomorphism $g= f\circ h_1\circ h_2\circ \cdots \circ h_\ell $ belongs to $D_{\alpha,\beta}$,
and we have: $$\De_{g,\varphi}=\psi_\ell\circ\cdots\circ\psi_1 \circ\De_{f,\varphi}.$$ 
\end{coro}

\begin{defi}\label{d.Theta}
Let $a\in\RR$, and let $\bar a$ be its projection on $S^1=\RR/\ZZ$. Given a diffeomorphism $\psi\colon S^1\to S^1$
with support in $S^1\setminus\{\bar a\}$ we call \emph{the lift of $\psi$  in $(a,a+1)$}
the diffeomorphism $\tilde \psi_a\colon\RR\to \RR$ with support in $(a,a+1)$ such that for any
$x\in(a,a+1)$ the image $\psi_a(x)$ is the point of $(a,a+1)$ which projects to $\psi(\bar x)$
where $\bar x$ is the projection of $x$. 

We denote by $\Theta_a(\psi)$ the diffeomorphism of $[0,1]$ whose expression in $(0,1)$
is $\Theta_a(\psi)=\varphi^{-1}\circ \psi_a\circ\varphi$.
\end{defi}

\begin{lemm}\label{l.voisinages} For any $C^1$ neighborhood $\cU$ of $f$ there is a neighborhood $\cV$
of $Id_{S^1}\in\Diff^1(S^1)$ with the following property:

Given any finite sequence $a_i>K_f$, $i\in\{1,\dots,\ell\}$, such that $a_{i+1}>a_i+1$
for all $i\in\{1,\dots,\ell-1\}$, we denote by $\bar a_i$ the projection of $a_i$ on $S^1$.
For any~$i$, let $\psi_i\in\cV$ be a diffeomorphism of $S^1$ with support in $S^1\setminus\{a_i\}$.
Then the diffeomorphism $g= f\circ \Theta_{a_1}(\psi_1)\circ \cdots \Theta_{a_\ell}(\psi_\ell)$
belongs to $\cU$.
\end{lemm}
\begin{demo} We fix a neighborhood $\cU_0$ if the identity map of $[0,1]$ such that,
if $g_1,\dots,g_n\in\cU_0$ and if the support of the $g_i$ are pairwise disjoint, then 
$ f\circ g_1\circ h_2\circ \cdots g_n $ belongs to $\cU$.  Now the lemma is a direct consequence
of Lemma~\ref{l.voisinages2} below.
\end{demo}
\begin{lemm}\label{l.voisinages2}
For any $C^1$-neighborhood $\cU_0$ of $f$ there is a neighborhood $\cV$ of $Id_{S^1}\in\Diff^1(S^1)$
with the following property:

Consider any  $a>K_f$, its the projection $\bar a$ on $S^1$ and any diffeomorphism  $\psi\in\cV$
with support in $S^1\setminus\{a\}$. Then the diffeomorphism $\Theta_{a}(\psi)$ belongs to $\cU_0$.
\end{lemm}
\begin{demo} Notice that there exists $\varepsilon>0$ such  that $\cU_0$ contains any diffeomorphism $h$ of $[0,1]$
with $\sup_{x\in[0,1]}|D_xh -1| <\varepsilon$. 
 
Now consider $a>K_f$ and an integer $n>0$. Then for any diffeomorphism $\psi$ of $S^1$ with support
in $S^1\setminus \{\bar a\}$, the lifts $\psi_{a}$ and $\psi_{a+n}$ are conjugated by the translation
$x\mapsto x+n$. As a consequence, $\Theta_{a+n}(\psi)$ is obtained from $\Theta_{a}(\psi)$ by the conjugacy
by the homothety of ratio  $\beta^n$.
As a consequence one gets that
$\sup_{x\in[0,1]}|D_x\Theta_{a+n}(\psi)-1|=\sup_{x\in[0,1]}|D_x\Theta_{a}(\psi)
-1|$.  

Hence one just has to prove the lemma for $a\in[K_f,K_f+1]$. This is a direct consequece of the facts
that the derivatives of $\varphi$ and of $\varphi^{-1}$ are bounded on 
$\varphi^{-1}([K_f,K_f+2])$ and $[K_f,K_f+2]$ respectively, and that for any $\psi$ with support in
$S^1\setminus\{a\}$, one has:
$$\sup_{x\in[a,a+1]}|D_x\psi_a-1|=\sup_{x\in S^1}|D_x\psi-1|.$$
\end{demo}

Let us now recall a classical result which is the key point of our proof.

\begin{theo}\label{t.decomposition}
Let $M$ be a closed Riemannian manifold, let $r>0$ and let $\cU$ be a $C^1$ neighborhood of the identity map.
Then for any smooth diffeomorphism $f$ of $M$ isotopic to the identity,
there exist $k\geq 1$ and $g_1,\dots , g_k\in \cU$ such that $g_i=id$ 
 on the complement of a ball $B(x_i,r)$, and 
$$f=g_1\circ\cdots\circ g_k.$$
\end{theo}
Here we use Theorem~\ref{t.decomposition} on the circle $S^1$, where it is an easy consequence of the result, 
by M. Herman, that any smooth diffeomorphism is the product of a rotation by a diffeomorphism 
smoothly conjugate to a rotation.
In Section~\ref{ss.sphere}, we will also use Theorem~\ref{t.decomposition} on the torus $T^2$.
\bigskip

\begin{demo}[\noindent Proof of Proposition~\ref{p.contre-exemple}]
Given a $C^1$-neighborhood $\cU$ of $f$, we choose a $C^1$-neighborhood $\cV$ of the identity map of $S^1$
given by Lemma~\ref{l.voisinages}. Using
Theorem~\ref{t.decomposition},
we can write $\De_{f,\varphi}$ as a finite
product 
$\De_{f,\varphi}=\psi_1^{-1}\circ\cdots\circ\psi_\ell^{-1}$ such that 
$\psi_i\in\cV$, and the support of $\psi_i$
is contained in an interval of length $\frac12$ in $S^1$ (and in
particular is not all of $S^1$).
Now we choose a finite sequence $a_i>K_k$ such that $a_{i+1}>a_i+1$, and such that the projection $\bar a_i$
does not belong to the support of $\psi_i$. Let $h_i=\Theta_{a_i}(\psi_i)$. 

Applying Lemma~\ref{l.voisinages}, we obtain that the diffeomorphism
$$g=f\circ h_1\circ h_2\circ \cdots \circ h_\ell $$
belongs to $\cU$; applying Corollary~\ref{c.composition}, we get that
$$\De_{g,\varphi}=\psi_\ell\circ\cdots\circ\psi_1 \circ\De_{f,\varphi}= Id_{S^1}.$$
\end{demo}

\subsubsection{Proof of Theorem~\ref{t=main2} on the circle}\label{sss=circle}
By Lemma~\ref{l.reduction}, it is enough to consider $f\in \cD^1_0$.
The set $Per(f)$ is finite. Let $\cI$ be the set of segments joining two successive periodic points of $f$;
in other words, every element $I\in\cI$ is the closure of a connected component of $S^1\setminus \hbox{Per}(f)$.
Notice that $f$ induces a permutation on $\cI$. Furthermore, all the elements of $\cI$ have the same period
denoted by $k>0$, under this action (this period is equal to $2$ if $f$ reverses the orientation,
and is equal to the period of the periodic orbits in the orientation preserving case).

Consider a segment $I\in \cI$. The endpoints of $I$ are the fixed points of the restriction $f^k|_I$;
moreover, one endpoint (denoted by $a$) is a repeller and the other (denoted by $b$) is an attractor.
Let $h_I\colon I\to [0,1]$ be the affine map such that $h_I(a)= 0$ and $h_I(b)=1$ and
let $\varphi_I\colon[0,1]\to [0,1]$ denote the diffeomorphism $h_I\circ f^k|_I\circ h_I^{-1}$.

According to Proposition~\ref{p.contre-exemple}, there is a sequence $(\psi_{I,n})_{n\in\NN}$, of diffeomorphisms
converging to $\varphi_I$ in the $C^1$-topology when $n\to+\infty$,
and a sequence $(Y_{I,n})_{n\in \NN}$ of $C^\infty$ vector fields on $[0,1]$
such that $\psi_{I,n}$ coincides with $\varphi_I$ in a small neighborhood of $\{0,1\}$
and is time one map of $Y_{I,n}$.
One denotes $g_{I,n}=h_I^{-1}\circ \psi_{I,n}\circ h_I$.
Notice that the diffeomorphism $g_{I,n}$ coincides with $f^k$ in 
neighborhoods of the endpoints of $I$ and
converges to $f^k|_I$ when $n\to \infty$.  

We now define a diffeomorphism $f_{I,n}$ of $S^1$ as follows:
$$f_{I,n}= 
\begin{cases}
f & \hbox{ on } S^1\setminus f^{k-1}(I)\\
g_{I,n}\circ f^{-k+1}&\hbox{ on }f^{k-1}(I).
\end{cases} $$
This is a $C^\infty$ diffeomorphism since it coincides with $f$ in a neighborhood of the
periodic orbits. Moreover, $(f_{I,n})$ converges to  $f$ when $n$ goes to $+\infty$.

One denotes by $X_{I,n}$ the vector field, defined on the orbit $\bigcup_0^{k-1}f^i(I)$
of the segment $I$ as follows:
\begin{itemize}
\item $X_{I,n}= (h_I^{-1})_*(Y_{I,n})$ on $I$;
\item for all $i\in\{1,\dots,k-1\}$ and all $x\in I$:
$$X_{I,n}(f^i(x))=f^i_*(X_{I,n}(x)).$$
\end{itemize}

Finally, we fix a family $I_1,\dots,I_\ell\subset \cI$ such that for $i\neq j$
the segments $I_i$ and $I_j$ have distinct orbits, and  
conversely every orbit of segment in $\cI$ contains one of the $I_i$. 

We denote by $f_n$ the diffeomorphism of $S^1$ coinciding with $f_{I_i,n}$ on the orbit of $I_i$
for all $i\in\{1,\dots, \ell\}$. This diffeomorphism is well-defined because all the  $f_{I_i,n}$
coincide with $f$ in a small neighborhood of the periodic points (the endpoints of the segments in $\cI$).
We denote by $X_n$ the vector field on $S^1$ that coincides with $X_{I_i,n}$ on the orbit of $I_i$,
for all $i\in\{1,\dots, \ell\}$.

One easily verifies that $X_n$ is a smooth vector field on $S^1$, invariant by $f_n$, and
such that $f_n^k$ is the time one map of $X_n$: the unique difficulty
consists in checking the continuity
and smoothness of the vector field $X$ at the periodic points.
As $f_n$ is affine in the neighborhood
of the periodic orbits one verifies that, at both sides of a periodic point $x$, the vector field $X$
is the affine vector field vanishing at $x$ and whose eigenvalue at
$x$ is  $\ln D_xf$. 

Finally $f_n$ converges to $f$ in the $C^1$ topology, completing the proof of Theorem~\ref{t=main2} on the circle.

\subsection{The case of the sphere $S^2$}\label{ss.sphere}
As in the one-dimensional case, the idea here is to
measure how far certain diffeomorphisms of $S^2$ are from the time-one map of
a vector field. One obtains in this way a 
generalization of the Mather invariant, which in this setting
is a diffeomorphism of $\TT^2$.
Such an invariant has already been constructed
\footnote{In~\cite{AY}, the authors
  write that the Mather invariant for a diffeomorphism of $S^2$ is always
isotopic to the identity, but this is not correct (their Proposition 1
is wrong). For this reason, we choose here to
build in detail the construction of this invariant on the sphere.} in~\cite{AY} by V. Afraimovich and T. Young,
and we now have to show that by a $C^1$-small perturbation of the dynamics, this invariant
vanishes.

\subsubsection{Preparation of diffeomorphisms in $\cO$}
Let $S^2$ be the unit sphere in $\RR^3$ endowed with the coordinates $(x,y,z)$. We denote by $N=(0,0,1)$ and $S=(0,0,-1)$ the north and the south poles of $S^2$. 
Notice that the coordinates $x,y$ define local coordinates of $S^2$ in local charts $U_N$ and $U_S$ in neigborhoods of $N$ and $S$

The following straightforward lemma asserts that one may assume that the fixed points of any diffeomorphism
$f$ in the open set $\cO$ are $N$ and $S$ and that the derivative at these points are confomal maps.  

\begin{lemm}\label{l.preparation}
Consider a diffeomorphism $f\in\cO$. Then there is a smooth
diffeomorphism $h\colon S^2\to S^2$ such that $h(N_f)=N$, $h(S_f)=S$
are the fixed points of $g=hfh^{-1}$; furthermore, the derivatives
$D_Ng$ and $D_Sg$ are conformal linear maps, i.e., each a composition
of a  rotation with  a homothety of ratio $\alpha>1$ and $\beta<1$, respectively.  

Finally, any $C^1$ neighborhood of $g$ contains a diffeomorphism $\tilde g$ such that there are neighborhoods $V_N\subset U_N$ and $V_S\subset U_S$ of $N$ and $S$, respectively, such that the expression of $\tilde g$ in the coordinates $(x,y)$ is $\tilde g(x,y)=D_Ng(x,y)$ for $(x,y)\in V_N$ and $\tilde g(x,y)=D_Sg(x,y)$ for $(x,y)\in V_S$.
\end{lemm}

\subsubsection{Space of orbits of a conformal linear map}

Let $A,B\in GL(\RR,2)$ be two conformal matrices of norm $\alpha>1$ and $\beta<1$, respectively.  There exist
$a,b\in[0,2\pi)$ such that $A=R_a\circ h_\alpha$ and $B=R_b\circ h_\beta$, where $R_a$ and $R_b$
are the rotation of angles $a$ and $b$, respectively, and $h_\alpha$
and $h_\beta$ are the homotheties of ratio $\alpha$ and $\beta$,
respectively. Notice that, for all $n\in\ZZ$, the linear map $A$ is
the time one map of the vector field 
$$X_{A,n}=\ln\alpha.(x\frac\partial{\partial
  x}+y\frac\partial{\partial y})+ (a+2\pi n).(x\frac\partial{\partial
  y}-y\frac\partial{\partial x}),$$ and $B$ is the time one map of 
$$X_{B,n}=\ln\beta.(x\frac\partial{\partial x}+y\frac\partial{\partial y})+ (b+2\pi n).(x\frac\partial{\partial y}-y\frac\partial{\partial x}).$$

Notice that the orbit space $T_A = \RR^2\setminus\{0\}/A$ (of the
action of $A$ on $\RR^2\setminus\{0\}$) is a torus (diffeomorphic to
$T^2=\RR^2/\ZZ^2$); we denote by $\pi_A$ the canonical projection from
$\RR^2\setminus\{0\}$ onto $T_A$. Moreover, the vector fields 
$$Z= 2\pi(x\frac\partial{\partial y}-y\frac\partial{\partial x})$$ and
$X_{A,n}$ project on $T_A$ in pairwise transverse commuting vector
fields, which we also denote by $Z$ and $X_{A,n}$; the orbits of both
flows are periodic of period $1$. Hence, for any pair $(Z,X_{A,n})$
there is  a diffeomorphism $\cL_{A,n}\colon T_A\to T^2=\RR^2/\ZZ^2$
which sends $Z$ to 
$\partial/\partial x$ and $X_{A,n}$ to $\partial/\partial y$; this
diffeomorphism is unique up to composition by a translation of $T^2$. 
Furthermore, the diffeomorphisms $\cL_{A,m}\circ \cL_{A,n}^{-1}$ are affine maps of the torus $T^2$, for all $n,m\in\ZZ$, so that $T_A$ is endowed with a canonical affine structure (indeed the affine map $\cL_{A,m}\circ \cL_{A,n}^{-1}$ is isotopic to the map induced by the matrix $\left(\begin{array}{cc} 1&n-m\\0&1\end{array}\right)$). 

Note that the orbits of $Z$ correspond to the positive generator of
 the fundamental group of $\RR^2\setminus \{0\}$; we denote by
 $\sigma$ the corresponding element of $\pi_1(T_A)$. Given any closed
 loop $\gamma\colon[0,1]\to T_A$, and any point
 $x\in\RR^2\setminus\{0\}$ with $\pi_A(x)=\gamma(0)$,
 there is a lift of $\gamma$ to a path in $\RR^2\setminus\{0\}$ 
joining $x$ to $A^k(x)$, 
where $k$ is the algebraic intersection number of $\sigma$ with $\gamma$. 
Finally, observe that the homotopy classes corresponding to the orbits of $X_{n,A}$, when $n\in\ZZ$ are precisely those whose intersection number with $\sigma$ is $1$: in other words, there is a  basis of $\pi_1(T_A)=\ZZ^2$ such that $\sigma=(1,0)$ and the orbits of $X_{A,n}$ are homotopic to $(n,1)$.

In the same way $T_B = \RR^2\setminus\{0\}/B$  is a torus endowed with the vector fields obtained by projection of $Z$ and $X_{B,n}$ and we denote by $\cL_{B,n}\colon T_B\to T^2$ a diffeomorphism that sends $Z$ to $\partial/\partial x$ and $X_{B,n}$ to $\partial/\partial y$.

\subsubsection{Mather invariant for diffeomorphisms of $S^2$}

Denote by $D_{A,B}\subset \cO$ the set of diffeomorphisms $f\in \cO$ whose expression in the coordinates $(x,y)$ coincides with $A$ in a neighborhood $U^N_f$ of $N$ and with $B$ in a neighborhood $U^S_f$ of $S$.
The aim of this part is to build a Mather invariant for diffeomorphisms in $D_{A,B}$.

We retain the notation of the previous subsection. 
Consider $f\in D_{A,B}$. The orbit space $\left(S^2\setminus\{N,S\}\right)/ f$ is a torus $T_f$ and we denote by $\pi_f\colon S^2\setminus\{N,S\}\to T_f$ the natural projection. Furthermore, as $f$ coincides with $A$ on $U^N_f$, the torus $T_f$ may be identified with the torus $T_A$ by a diffeomorphism $\pi_N\colon T_f\to T_A$, and in the same way, the fact that $f$ coincides with $B$ in a neighborhood of $S$ induces a diffeomorphism $\pi_S\colon T_f\to T_B$. 

Notice that the morphisms $\pi_{N*} \colon H_1(T_f,\ZZ) \to H_1(T_A,\ZZ)$
and $\pi_{S*} \colon$ $H_1(T_f) \to H_1(T_B)$ preserve the homology class of $\sigma$ (corresponding to the positive homology generator of $S^2\setminus \{N,S\}$ or of $\RR^2\setminus \{0\}$), and the homology intersection form with $\sigma$. 

Consequently, for any $f\in D_{A,B}$, there is an integer $n(f)$ such that  the map $\De_{f,0,0}=\cL_{B,0}\circ \pi_S\circ \pi_N^{-1} \circ\cL_{A,0}^{-1}$ is isotopic to the linear map of $T^2$ induced by the matrix $\left(\begin{array}{cc} 1&n(f)\\0&1\end{array}\right)$.

\begin{lemm} For any $f\in D_{A,B}$ there is a $C^1$-neighborhood $\cU$ of $f$ in $\Diff^1(S^2)$ such that for any $g\in \cU\cap D_{A,B}$ one has $n(f)=n(g)$. 
\end{lemm}
\begin{demo} 
We can choose a neighborhood $\cU$ such that, if $g\in\cU$ then the map
$$f_t(x)= \frac{(1-t)f(x)+t g(x)}{\|(1-t)f(x)+t g(x)\|}$$ 
is a smooth isotopy between $f$ and $g$. 
Furthermore, by shrinking $\cU$ if necessary, for any 
$g\in\cU$, the isotopy $f_t$ belongs to $\cO$ (that is $\Om(g)=\{N_g,S_g\}$). 

If $g\in \cU\cap D_{A,B}$ then there are discs $D^N$ and $D^S$ centered on $N$ and $S$, respectively, such that $f_t=A$ on $D^N$ and $f_t=B$ on $D^S$ so that $f_t\in D_{A,B}$. In particular $f_t(D^S)\subset D^S$, and $f_t^{-1}(D^N)\subset D^N$. Furthermore, there exists $\ell>0$ such that for any $x\in S^2\setminus (D^N\cup D^S)$,  $f_t^\ell(x)\in D^S$ and $f_t^{-\ell}(x)\in D^N$. 

Let $x\in D^N$ such that $A(x)=f_t(x)\in D^N$ and $A^2(x)\notin
 D^N$. Hence $y_t=f_t^{\ell+2}(x)\in D^S$ and $f_t(y_t)=B(y_t)\in
 D^S$. Let $\gamma$ be the segment of orbit of $X_{A,0}$  joining $x$
 to $A(x)=f_t(x)$, and let $\gamma_t= f_t^{\ell+2}(\gamma)$.
For every $t$, $\gamma_t$ is homotopic (relative to $\{y_t,B(y_t)\}$
in $S^2\setminus \{N,S\}$) to a segment of orbit of $X_{B,n(f_t)}$. As a consequence, $n(f_t)$ varies continuously with $t$ as $t$ varies from $0$ to $1$. Hence $n(f_t)$ is constant; that is, $n(g)=n(f)$.  
\end{demo}

Hence there is a partition of $D_{A,B}$ into open subsets $D_{A,B,n}$ 
such that $n(f)=n$ for $f\in D_{A,B,n}$. For $f\in D_{A,B,n}$, we
define:
$$\De_f=\cL_{B,n}\circ \pi_S\circ \pi_N^{-1} \circ\cL_{A,0}.$$ 
Then $\De_f$ is a diffeomorphism of $T^2$, isotopic to the identity. 

Theorem~\ref{t.mather2dim} below justifies calling $\De_f$ \emph{the Mather invariant of $f$}. 

\begin{theo}\label{t.mather2dim} Let $f\in D_{A,B,n}$ be a smooth diffeomorphism such that $\De_f$ is a translation of the torus $T^2$. Then $f$ leaves invariant two transverse  commuting vector fields $Z_f$ and $X_f$ on $S^2$ such that $Z_f=Z$ in a neighborhood of $\{N,S\}$, $X_f=X_{A,S}$ in a neighborhood of $N$ and $X_f=X_{B,S}$ in a neighborhood of $S$. 

As a consequence the centralizer of $f$ is isomorphic to $S^1\times \RR$. 
\end{theo} 
\begin{demo} Fix two discs  $D^N$ and $D^S$ centered at  $N$ and $S$,
  respectively,  in which $f$ coincides with $A$ and $B$, respectively. 

For any $x\neq S$ there exists $m(x)<0$ such that $f^{m(x)}\in
D^N$. One defines $Z_f(x)=f^{-m(x)}_*(Z(f^{m(x)}))$ and
$X_f(x)=f^{-m(x)}_*(X_{A,N}(f^{m(x)}))$. As $Z$ and $X_{A,N}$ are
invariant by $A$, one proves that the vectors $Z_f(x)$ and $X_f(x)$
are independent of the choice of $m(x)$. As a consequence, one deduces
that they depend smoothly on $x\in S^2\setminus \{S\}$ and that they 
commute on $S^2\setminus\{S\}$. Furthermore the restrictions of $Z_f$
and $X_f$ to $D^S$ are invariant by $f$, and hence by $B$, so that
they induce two vector fields on $T_B$ whose images by $\cL_{B,n}$ are
$\De_f(\frac{\partial}{\partial x})= \frac{\partial}{\partial x}$ and
$\De_f(\frac{\partial}{\partial y})=\frac{\partial}{\partial y}$,
respectively; that is, they
coincide with the projections of the restrictions $Z$ and $X_{B,S}$ to
$D^S$. 
Thus $Z_f=Z$ and $X_f=X_{B,S}$ on $D^S$, ending the proof. 
\end{demo}

\subsection{Vanishing of the Mather invariant}
This part is now very close to the $1$-dimensional case. 

For any $f\in D_{A,B,n}$ we denote by $D^S_f$ a disk centered on $S$ on which $f=B$. 

Let $h\colon S^2\to S^2 $ be a diffeomorphism whose support is contained in a disk $D\subset D^S_f$, disjoint from all $B^m(D)$ for $m>0$. The disk $D$ projects homeomorphically onto a disk  $D'\subset T_B$, and finally onto a disk $\tilde D=\cL_{B,n}(D')\subset T^2$. Let $\psi$ be the diffeomorphism of $T^2$ with support in $\tilde D$ whose restriction to $\tilde D$ is the projection of $h$.  We says that $\psi$ is the projection of $h$ on $T^2$ and conversely, that $h$ is the lift of $\psi$ with support in $D$. 

Fix $k>0$ such that $D$ is disjoint from $B^k(D^S_f)$.

\begin{lemm} With the notation above, the composition $f\circ h$ is a diffeomorphism in $D_{A,B,n}$ with $B^k(D^S_f)\subset D^S_{f\circ h}$, and whose Mather invariant is $$\De_{f\circ h}= \psi\circ \De_f.$$ 
\end{lemm}     
\begin{coro} Let $D_0,\dots, D_\ell\subset D^S_f$ be a finite sequence of disks such that 
\begin{itemize} 
\item for every $i$, $D_i$ is disjoint from $B^k(D_i)$ for $k>0$;
\item for all $i<j$ the disk $D_i$ is disjoint from $B^k(D_j)$, $k\geq 0$
\end{itemize}

For every $i$, let $h_i$ be a diffeomorphism of $S^2$ with support in
$D_i$, and let $\psi_i$ be the projection of $h_i$ on $T^2$ (by $\cL_{B,n}\circ\pi_S\circ\pi_f$).

Then the Mather invariant of $f\circ h_0\circ\dots\circ h_\ell$ is 
$$\De_{f\circ h_0\circ\dots\circ h_\ell}= \psi_\ell\circ\dots\circ \psi_0\circ \De_f.$$ 
\end{coro} 

Observe that for any disk $\tilde D\subset T^2$ with diameter strictly
less than $1$, each connected component of 
$(\cL_{B,n}\circ\pi_S\circ\pi_f)^{-1}(\tilde D)$ projects 
diffeomorphically onto $\tilde D$, and $f$ induces a permutation of
these components. For $i>0$, let $D_i$ denote the (unique) component
of $(\cL_{B,n}\circ\pi_S\circ\pi_f)^{-1}(\tilde D)$ such that $
f^{-i}(D_i)\subset D^S_f$ but $f^{-(i+1)}(D_i)$ is not contained in $D^S_f$.

For any diffeomorphism $\psi$ with support in $\tilde D$  we will
 denote
by  $\theta_i(\psi)\colon S^2\to S^2$ the lift of $\psi$ with support in $D_i$. 

The next lemma is the unique reason we required that the derivative of $f$ at $N,S$ be complex, hence conjugate to conformal linear maps:

\begin{lemma} Let $\tilde D\subset T^2$ be a disk with diameter
  strictly less than $1$ and let $i,j\in\NN$. Then:
$$ \sup_{x\in S^2 }\|D_x\theta_i(\psi)-Id\| =\sup_{x\in S^2}\|D_x\theta_j(\psi)-Id\|.$$
\end{lemma}
\begin{demo} $\theta_i(\psi)$ is conjugated to $\theta_j(\psi)$ by $B^{j-i}$ which is the composition of a homothety by a rotation; the $C^1$ norm is preserved by conjugacy by an isometry, and also by conjugacy by a homothety, hence is preserved by the conjugacy by $B^{j-i}$. 
\end{demo}

\begin{coro}\label{c.varep} For any $\varepsilon>0$ there is a $C^1$-neighborhood $\cV_\varepsilon\subset \Diff(T^2)$ of the identity map such that for any diffeomorphism $\psi\in\cV_\varepsilon$ with support in a disk $\tilde D\subset T^2$ with diameter strictly less than $1$, and for any $i\geq 0$, the lift $\theta_i(\psi)$ satisfies : 
$$ \sup_{x\in S^2}\|D_x\theta_i(\psi)-Id\| <\varepsilon.$$
\end{coro}

\begin{defi} Let $\psi_1,\dots,\psi_\ell$ be $\ell$ diffeomorphisms of
$T^2$ such that the support of every $\psi_i$ is contained in a disk
$\tilde D_i$ with diameter strictly less than $1$;  a {\em lift of the sequence $\psi_1,\dots,\psi_\ell$} is a sequence of lifts $h_1= \theta_{i_1}(\psi_1),\dots,h_\ell= \theta_{i_\ell}(\psi_\ell)$ such that, for every $i<j$ the support of $h_i$ is disjoint from all the iterates $B^k(supp(h_j))$, for $k\geq 0$. 
\end{defi}

It is easy to check that, for any sequence $\psi_1,\dots,\psi_\ell$ of diffeomorphisms of $T^2$ such that the support of every $\psi_i$ is contained in a disk $\tilde D_i$ with diameter strictly less than $1$, the sequence $h_i=\theta_{i}(\psi_i)$ is a lift.
\bigskip

\begin{demo}[\noindent Proof of Theorem~\ref{t=main2} on the sphere $S^2$] 
Consider  $f\in D_{A,B,n}$ and a $C^1$-neighborhood $\cU$ of $f$. Fix
$\varepsilon>0$ such that, if $g_1,\dots,g_m$, $m>0$, are
diffeomorphisms of $S^2$ with pairwise disjoint supports in
$S^2\setminus \{N,S\}$, and such that   
$ \sup_{x\in S^2} \|Dg_i(x)-Id\| <\varepsilon$, then $f\circ
g_1\circ\cdots\circ g_m\in\cU$. Let $\cV_\varepsilon$ be the
$C^1$-neighborhood
of the identity map of $T^2$ given by Corollary~\ref{c.varep}.

Using Theorem~\ref{t.decomposition}, we write 
$$\De_f=\psi_1^{-1}\circ\cdots\circ \psi_\ell^{-1},$$
for some $\ell>0$, where $\psi_i\in \cV_\varepsilon$, and 
the support of $\psi_i$ is contained in a disk $\tilde D_i$ 
with diameter strictly less
than $1$. 
Let $(h_1,\dots,h_\ell)$ be a lift of the sequence 
$(\psi_1,\dots,\psi_\ell)$; the $h_i$ satisfy 
$$ \sup_{x\in S^2} \|D_xh_i-Id\| <\varepsilon,$$ 
by our choice of $\cV_\varepsilon$. 

Our choice $\varepsilon>0$ implies that 
$g=f\circ h_1\circ\cdots\circ h_\ell$ is a diffeomorphism belonging
to $D_{A,B,n}\cap\cU$. Furthermore, its Mather invariant is 
$\De_g=\psi_\ell\circ\cdots\circ\psi_1\circ \De_f=Id$. 

We have just shown that any $f\in D_{A,B}$ is the $C^1$-limit of a
sequence $g_k\in D_{A,B}$ whose Mather invariant is the identity map;
in particular, the centralizer of $g_k$ is isomorphic to $\RR\times S^1$. 

Since by Lemma~\ref{l.preparation}, $\cO$ contains a dense set of diffeomorphisms smoothy conjugate to elements of $D_{A,B,n}$, any diffeomorphism in $\cO$ is the limit of diffeomorphisms $g_k$ that are the time $1$ map
of Morse-Smale vector fields, ending the proof of Theorem~\ref{t=main2}. 
\end{demo}

\appendix

\section{Conservative extension results}\label{a=invmanif}
We explain in this appendix how a perturbation of a conservative diffeomorphism
along a submanifold $W$ can be extend as a conservative perturbation on the whole manifold $M$.

This allows to obtain Proposition~\ref{p=main} from Proposition~\ref{p=diskdist}:
the results proven in this section will be applied to the case $W$ is an invariant manifold
of a hyperbolic periodic point $p$. In the volume-preserving case, one will
assume that $\dim(W)\leq \frac 1 2 \dim(M)$ (note that this hypothesis is always satisfied either
by the stable or by the unstable manifold of $p$).
In the symplectic case, there is no additional hypothesis, but we use the following well-known fact.

\begin{lemma} Let $f\in \hbox{Symp}^1(M)$ and let $p$ be a hyperbolic periodic point for $f$.
Then $W^s(p)$ and $W^u(p)$ are Lagrangian submanifolds of $M$.
\end{lemma}
\begin{proof}
Let $x\in W^s(p)$, and let $v,w\in T_xW^s(p)$ be
tangent vectors to $W^s(p)$. On the one hand, since $f$ is a symplectomorphism,
we have
$$\omega(D_xf^k(v), D_xf^k(w)) = \omega(v,w),$$
for all $k\in \ZZ$.
On the other hand, as $k\to +\infty$, we have
$$\omega(D_xf^k(v), D_xf^k(w))\to 0.$$
Hence $\omega$ vanishes identically on $W^s(p)$.
The same is true for $W^u(p)$. Since $W^s(p)$ and $W^u(p)$ have
complementary dimension and $\omega$ is nondegenerate, they must have the same dimension.
Hence, both are Lagrangian submanifolds of $M$.
\end{proof}

\subsection{The symplectic case}\label{s=symplectic}
\begin{prop}\label{p.sympert} 
Let $M$ be a symplectic manifold and $z$ a point contained in a $C^1$ Lagrangian
submanifold $W\subset M$. Then there exists in $W$ a disk $D=\overline{B_W(z,r_0)}$ centered at $x$ such that,
for every neighborhood $U\subset M$ of $D$ and every $\eps>0$,
there exists $\delta>0$ with the following property.

For every $C^1$ diffeomorphism $\psi:D\to D$ satisfying:
\begin{itemize}
\item[a.] $\psi=Id$ on a neighborhood of $\partial D$, and  
\item[b.] $d_{C^1}(\psi, Id) <\delta$,
\end{itemize}
there exists $\varphi\in\hbox{Symp}^1(M)$ such that:
\begin{enumerate}
\item $\varphi=Id$ on $M\setminus U$,
\item $\varphi = \psi$ on $D$, and
\item $d_{C^1}(\varphi,Id) <\eps$.
\end{enumerate}
\end{prop}
\begin{proof}
The basic strategy is first to symplectically
embed the disk $D$ as the $0$-section of its cotangent bundle
$T^*D$. On $T^*D$, the symplectic form is $\omega = d\alpha$, where $\alpha$ is the canonical one-form on $T^*D$.  Any diffeomorphism $\psi:D\to D$ 
lifts to a canonical symplectomorphism $\psi^\ast: T^\ast D\to T^\ast D$;
namely the pull-back map $(\psi, D\psi^{-1})$.  The natural thing to
try to do is to set 
$\varphi = \psi^\ast$ in a neighborhood of the $0$-section, symplectically
interpolating between $\psi^\ast$ and $Id$ using a generating function.
This simple approach fails, however, because $\psi$ is only $C^1$,
and so $\psi^\ast$ is merely continuous.
(Even assuming that $\psi$ is $C^2$ does not help:
in order to control the $C^1$ size
of such a map, it is necessary to have some control on the $C^2$ size of
$\psi$, and we cannot assume any such control).  Using a convolution product, 
it is possible to overcome this problem.   This approach mirrors that in 
\cite{BGV}, but in the symplectic setting.

The problem is local and one can work in $\RR^{2n}$ endowed with the standard symplectic form
$\omega = \sum_i du_i\wedge dv_i$ where $u = (u_1,\ldots, u_n), v=(v_1,\ldots, v_n)$.
By a symplectic change of coordinates, we may assume that the disk $D$ lies inside
a disk $\{(u,v),\,\|u\|\leq R, \, v=0\}$.
We define $\psi$ using a generating function $S$.  

We first recall the definition and properties of generating functions.
Suppose that $h:{\RR^{2n}}\to \RR^{2n}$ is a $C^r$ symplectomorphism, taking the form:
$$h(u,v) = (\xi(u,v),\eta(u,v)),$$
with $\xi,\eta:{\RR}^{2n}\to \RR^n$ and
$h(0,0) = (0,0)$.
Let us assume that the partial derivative matrix $\frac{\partial}{\partial v}\eta(u,v)$ is
invertible (this is the case for instance if $h$ preserves $\RR^n \times \{0\}$).
We can solve for $\eta = \eta(u,v)$ to obtain new coordinates
$(u,\eta)$ on a small neighborhood of $(0,0)$ in  $\RR^{2n}$.
Since $h$ is symplectic, the $1$-form $\alpha = \sum_i v_i du_i + \xi_i d\eta_i$ is closed, and hence, exact.
Thus there exists a $C^{r+1}$ function $S = S(u, \eta)$, unique up to adding a constant, defined in
a neighborhood of $(0,0)$, such that $dS = \alpha$.
The function $S$ is called a {\em generating function for $h$}.

On the other hand, any $C^{r+1}$ function $S=S(u,\eta)$ satisfying the nondegeneracy condition that $\frac{\partial^2}{\partial u\partial \eta}S$ is everywhere nonsingular is the generating function of a $C^r$ symplectic diffeomorphism.
Solving for $\alpha$ in the
equation
$$dS = \frac{\partial S}{\partial u} du + \frac{\partial S}{\partial\eta} d\eta = \alpha = v du + \xi d\eta ,$$
we obtain the system:
$$\frac{\partial S}{\partial u} = v;\qquad \frac{\partial S}{\partial\eta} = \xi.$$
The nondegeneracy condition implies that this system can be solved
implicitly for a $C^r$ function $\eta = \eta(u,v)$.  
We then obtain a  $C^r$ symplectomorphism:
$$h(u,v) = \left(\frac{\partial S}{\partial\eta}(u,\eta(u,v)), \eta(u,v)\right),$$
and $S$ is a generating function for $h$.

It is easy to see that the generating function for the identity map is
$$S_0(u,\eta) = u\cdot \eta = \sum_{i=1}^n u_i\eta_i.$$

\begin{affi}
For every $\eps>0$, there exists $\delta>0$ such that, 
if $d_{C^2}(S,S_0) < \delta$ then $d_{C^1}(h,Id) <\eps$.
\end{affi}

\noindent{\bf Proof.} This follows from the implicit function theorem,
and the details are omitted.\eproof

Returning to the proof of Proposition~\ref{p.sympert}, assume  
that $\psi:D\to D$ is written in $u$-coordinates as
$$\psi(u_1,\ldots, u_n) = (\psi_1(u_1,\ldots,u_n),\ldots, \psi_n(u_1,\ldots, u_n)).$$
We may assume that the domain of $\psi$ has been extended to 
$\RR^n$.  
To prove Proposition~\ref{p.sympert}, it suffices to prove the following lemma.
\end{proof}

\begin{lemma}\label{l.sympert} Given a disk $D\subset \RR^n$,
and a neighborhood $U$ of $D\times \{0\}$
in $\RR^{2n}$,  there exists $C>0$ with the following property.

For every $C^1$ diffeomorphism $\psi:\RR^n\to \RR^n$, equal
to the identity on a neighborhood of $\partial D$, there is a $C^2$
function $S:\RR^{2n}\to \RR$ such that:
\begin{enumerate}
\item $d_{C^2}(S_0,S) \leq Cd_{C^1}(\psi, Id)$,
\item $S=S_0$ outside of $U$,
\item $\frac{\partial S}{\partial u}(u,0) = 0$ for all $u\in \RR$ and
\item $\frac{\partial S}{\partial \eta}(u,0) = \psi(u)$ for all $u\in D$.
\end{enumerate}
\end{lemma}
Note that condition 1. implies that $S$ is nondegenerate, provided that
$d_{C^1}(\psi, Id)$ is sufficiently small.

\bigskip

\noindent{\bf Proof of Lemma~\ref{l.sympert}.}
To illustrate the argument in a simple case, we
first prove the lemma for $n=1$.
The proof of the general case is very similar.
Let
$$a(u) = \psi'(u)-1.$$ 
Note that $a$ is a continuous map, $\|a\|_\infty\leq d_{C^1}(\psi, Id)$, and
$a(u) =0$ if $u\notin \hbox{int}(D)$. 
Let $\Phi:\RR\to [0,1]$ be a $C^\infty$ function satisfying:
\begin{itemize}
\item $\Phi(0) = 1$ and $\Phi = 0$ outside of $(-1,1)$,
\item $\Phi^{(k)}(0) = 0$, for all $k\geq 1$,
\item $\int_\RR \Phi(w)\, dw = 1$.
\end{itemize}
Fix a point $u_\ast\in \partial D$, so that $\psi(u_\ast) = u_\ast$.
For $(u,\eta)\in \RR^{2}$, $\eta\neq 0$, let:
$$Q(u,\eta) = \eta \int_{u_\ast}^u\int_\RR \Phi(w)\, a(x-w\eta)\,dw\,dx.$$
For $\eta\neq 0$, one can make the change of variables
$w' = x-w\eta$ and get
$$Q(u,\eta)= \sign(\eta)\int_{u_\ast}^u \int_\RR
\Phi\left(\frac{x-w'}{\eta}\right) a(w')\,dw'\,dx.$$
Let $\rho:\RR^2\to \RR$ be a $C^\infty$ bump function
identically equal to $1$ on a neighborhood of $D\times \{0\}$
and vanishing outside of $U$. Consider
$$S= S_0 + \rho\, Q.$$
Lemma~\ref{l.sympert} in the case $n=1$ is a direct consequence of:

\begin{clai}\label{l.conv} The map $Q\colon \RR^{2}\to \RR$ is $C^2$ and 
there is $C = C(U) >0$ such that:
\begin{enumerate}
\item $\|Q\,\vert_{\overline U}\|_{C^2} \leq C \|a\|_\infty$,
\item $\frac{\partial Q}{\partial u}(u,0) = 0$, for all $u\in \RR$, and
\item $\frac{\partial Q}{\partial \eta}(u,0) = 
\int_{u_\ast}^u a(x)\,dx = \psi(u)-u$, for all  $u\in \RR$.
\end{enumerate}
\end{clai}
\begin{proof}
We derive explicitly the formulas: 
\begin{eqnarray*}
\frac{\partial Q}{\partial u} 
&=& \eta \int_\RR \Phi(w) \, a(u-w\eta)\,dw\\
&=& \sign(\eta)\int_{\RR} \Phi\left(\frac{u-w'}{\eta}\right)\, a(w')\,dw',\\
\\
\frac{\partial Q}{\partial \eta}& =&
\frac{-\sign(\eta)}{\eta^2}\int_{u_\ast}^u \int_{\RR} \Phi'\left(\frac{x-w'}{\eta}\right)\,(x-w') \,a(w')\,dw'\,dx\\
&=& -\int_{u_\ast}^u \int_\RR \Phi'(w)\, w\, a(x-w\eta)\, dw\,dx\\
&=&-\int_\RR \Phi'(w)\, w\,\int_{u_\ast - w\eta}^{u-w\eta} a(x')\,dx'\,dw,
\end{eqnarray*}
\begin{eqnarray*}
\frac{\partial^2 Q}{\partial \eta \partial u}& =& 
-\int_\RR \Phi'(w)\, w\, a(u-w\eta)\, dw,\\
\\
\frac{\partial^2 Q}{\partial u^2} &=&
\frac{1}{|\eta|}\int_{\RR} \Phi'\left(\frac{u-w'}{\eta}\right) \,a(w')\,dw'\\
&=& \int_\RR \Phi'(w)\,a(u -w\eta)\, dw,\\
\\
\hbox{and finally:}&&\\
\\
\frac{\partial^2 Q}{\partial \eta^2}
&=& \int_\RR \Phi'(w) \,w^2 \,\left(a(u-w\eta) - a(u_\ast -w\eta)\right) \,dw.
\end{eqnarray*}
Properties 1. and 2. follow immediately 
from these formulas. To see 3., note that
\begin{eqnarray*}
\frac{\partial Q}{\partial \eta}\vert_{\eta = 0}
&=& -\left(\int_{u_\ast}^u a(x)\,dx\right)  \left(\int_\RR \Phi'(w) w\, dw\right)\\
&=& -\left(\int_{u_\ast}^u a(x)\,dx\right)  \left(-\int_\RR \Phi(w) dw\right)\\
&=&\int_{u_\ast}^u a(x)\,dx.
\end{eqnarray*}
\end{proof}

We now turn to the case $n\geq 1$ in Lemma~\ref{l.sympert}.
For $i=1,\ldots n$, let $\alpha_i$ be the continuous $1$-form defined by
$$\alpha_i = d(\psi_i - \pi_i),$$
where $\pi_i:\RR^n\to \RR$ is the projection onto the $i$th coordinate.
As above, fix a point $u_\ast\in \partial D$, so that $\psi(u_\ast) = u_\ast$.
Then we have the formula:
$$\psi_i(u_1,\ldots,u_n) - u_i =
\int_{u_\ast}^u \alpha_i,$$
where the right-hand side is a path integral evaluated
on any path from $u_\ast$ to $u = (u_1,\ldots, u_n)$.
Furthermore, we have $\|\alpha_i\|_\infty \leq d_{C^1}(\psi,Id)$, for all $i$.
When $n=1$, the $1$-form $\alpha_1$ is just
$\alpha_1 = a(u)\, du$, where $a(u) = \psi'(u)-1$, as above.

Let $\Phi_n:\RR^n\to\RR$ be an 
$n$-dimensional bell function:
$$\Phi_n(x_1,\ldots, x_n) = \Phi(x_1)\cdots \Phi(x_n).$$
For each $1$-form $\alpha$, and $t\in \RR$, we define
a new $1$-form $\alpha_i^{\star t}$ on $\RR^n$ by taking the convolution:
$$\alpha^{\star t}(u) = t\int_{\RR^n}  \Phi_n(w)\, \alpha(u-tw) \, dw.$$

We integrate along any path from $u_\ast$ to $u$ and set
$$Q(u,\eta) = \int_{u_\ast}^u \sum_{i=1}^n \alpha_i^{\star \eta_i}
=\sum_{i=1}^n \eta_i\int_{\RR^n}\Phi_n(w) \,\left(\int_{u_\ast}^u
\alpha_i(u-tw)\right)\,dw.$$
This is well-defined since
$\int_{u_\ast}^u\alpha_i(u-tw)$
is independent of choice of path.

Let $\rho_n:{\bf R}^{2n} \to [0,1]$ be a $C^\infty$ bump
function vanishing identically outside of $U$ and equal to $1$
on a neighborhood of $D$.
As before, the map $S = S_0 +\rho_n Q$ satisfies the conclusions of Lemma~\ref{l.sympert}
provided the following claim holds.

\begin{clai} The map $Q\colon\RR^{2n}\to \RR$ is $C^2$ and there is $C = C(U) >0$ such that:
\begin{enumerate}
\item $\|Q\,\vert_{\overline U}\|_{C^2} \leq C \max_i\|\alpha_i\|_\infty$,
\item $\frac{\partial Q}{\partial u}(u,0) = 0$, for all $u\in \RR^n$, and
\item $\frac{\partial Q}{\partial \eta_i}(u,0) = 
\int_{u_\ast}^u \alpha_i  = \psi_i(u)-u_i$, for all $1\leq i\leq n$ and  $u\in \RR$.
\end{enumerate}
\end{clai}
\begin{proof}
We repeat the calculations from the proof of Lemma~\ref{l.conv} in the general setting.
When $t\neq 0$, the change of variable $w'=u-tw$ gives
\begin{eqnarray*}
\alpha^{\star t} (u) &=& \sign(t)\int_{\RR^n} \Phi_n
\left(u-tw'\right)\,\alpha(w' )\,dw',\\
\frac{d}{dt} \alpha^{\star t} (u) &=& -\int_{\RR^n}
\left(d\Phi_n(w).w+(n-1)\Phi_n(w)\right)\,\alpha(u-tw)\,dw.
\end{eqnarray*}
One deduces:
\begin{eqnarray*}
\frac{\partial Q}{\partial u} &=&
\sum_{i=1}^n \eta_i \int_{\RR^n} \Phi_n(w)\, \alpha_i(u-tw)\,dw\\
&=& \sum_{i=1}^n \frac{\eta_i}{|\eta_i|^n}\int_{\RR^n}
\Phi_n\left(\frac{u-w'}{\eta_i}\right)\, \alpha_i(w')\,dw',\\
\\
\frac{\partial Q}{\partial \eta_i}& =&
 \int_{u_\ast}^u \frac{d}{d\eta_i} \alpha_i^{\star\eta_i}\\
&=&-\int_{\RR^n} \left(d\Phi_n(w).w+(n-1)\Phi_n(w)\right)\,
\left(\int_{x=u_\ast-\eta_iw}^{x=u-\eta_iw}\alpha_i(w)\right) \, dw, \\
\\
\frac{\partial^2 Q}{\partial u\, \partial \eta_i}
&=&-\int_{\RR^n} \left(d\Phi_n(w).w+(n-1)\Phi_n(w)\right)\,
\alpha_i(u-\eta_iw)\, dw,\\
\\
\frac{\partial^2 Q}{\partial u^2} &=&
\sum_{i=1}^n 
\int_{\RR^n} d \Phi_n(w)\,\alpha_i(u-\eta_iw)\, dw,\\
\end{eqnarray*}
and finally:
\begin{eqnarray*}
\frac{\partial^2 Q}{\partial \eta_i \partial \eta_j}
&=&  \delta_{i,j} \int_{\RR^n} \left(d\Phi_n(w).w+(n-1)\Phi_n(w)\right)\,
\left[\alpha_i(x-\eta_iw).w\right]_{x=u_\ast}^{x=u} \, dw.\\
\end{eqnarray*}
It is not difficult to verify that 1.--3. hold.
\end{proof}

The proof of Lemma~\ref{l.sympert} is now complete.\endproof

\subsection{The volume-preserving case}\label{s=volpres}

\begin{prop}\label{p.vpert} Let $M$ be a Riemannian manifold
endowed with a volume form $\mu$ and $W$ be a $C^1$ submanifold satisfying
$$\hbox{dim}(W) \leq \hbox{codim}(W).$$
Centered at any point $z\in W$, there exists a disk $D=\overline{B_W(z,r_0)}$
of $W$ such that, for every neighborhood $U\subset M$ containing $D$ and every $\eps>0$,
there exists $\delta>0$ with the following property.  

For every $C^1$ diffeomorphism $\psi:D\to D$ satisfying:
\begin{itemize}
\item[a.] $\psi=Id$ on a neighborhood of $\partial D$, and  
\item[b.] $d_{C^1}(\psi, Id) <\delta$,
\end{itemize}
there exists $\varphi\in\Diff^1_\mu(M)$ such that:
\begin{enumerate}
\item $\varphi=Id$ on $M\setminus U$
\item $\varphi = \psi$ on $D$, and 
\item $d_{C^1}(\varphi,Id) <\eps$.
\end{enumerate}
\end{prop}
\begin{proof} Let $n=\hbox{dim}(M)$.
By a local change of coordinates, we may assume that
$\mu$ is the standard volume form $dx_1\wedge \cdots\wedge dx_n$
on a neighborhood of the origin in $\RR^n$.
By composing these coordinates with an isometry of $\RR^n$,
we may further assume that $D$ is the graph of a $C^1$-function $h:\RR^k \to \RR^{n-k}$,
where $k\leq n/2$.
The final change of coordinates
$$(x_1,\ldots, x_n)\mapsto ((x_1,\ldots, x_k), (x_{k+1},\ldots, x_n) - h(x_1,\ldots, x_k))$$
preserves volume.  Applying this change of coordinates, we may assume
that $D$ lies in the coordinate plane
$\{x_{k+1}=x_{k+2}=\cdots=x_n=0\}\simeq \RR^k$. Now we apply the symplectic
pertubation result (Proposition~\ref{p.sympert}) inside the
space $\{x_{2k+1}=\cdots=x_n = 0\}\simeq \RR^{2k}$ to obtain a local
$C^1$ symplectomorphism $\varphi_0$ of $\{x_{2k+1}=\cdots=x_n= 0\}$ that agrees with 
$\psi$ on $D$. This symplectomorphism is 
$C^1$-isotopic to the
identity through symplectomorphisms $\{\varphi_t\}_{t\in[0,1]}$,
where $\varphi_1 = Id$ (to obtain this
isotopy, just choose a smooth
isotopy of the generating function for $\psi$
to the generating function for the identity).
 
Now we extend $\varphi_0$ to $\RR^n$
using this isotopy to obtain a 
locally-supported volume-preserving diffeomorphism
that agrees with $\psi$ on $D$.  More precisely, choose an
 appropriate
$C^\infty$ bump function $\rho:\RR^{n-2k}\to [0,1]$, and set
$$\varphi(x_1,\ldots, x_n) = (\varphi_{\rho(\|(x_{2k+1},\ldots,x_n)\|)}(x_1,\ldots, x_{2k}), x_{2k+1},\ldots, x_n).$$
This is the desired map $\varphi$.
\end{proof}

\vspace{10pt}

\noindent \textbf{Christian Bonatti (bonatti@u-bourgogne.fr)}\\
\noindent  CNRS - Institut de Math\'ematiques de Bourgogne, UMR 5584\\
\noindent  BP 47 870\\
\noindent  21078 Dijon Cedex, France\\
\vspace{10pt}

\noindent \textbf{Sylvain Crovisier (crovisie@math.univ-paris13.fr)}\\
\noindent CNRS - Laboratoire Analyse, G\'eom\'etrie et Applications, UMR 7539,\\
\noindent Institut Galil\'ee, Universit\'e Paris 13, Avenue J.-B. Cl\'ement,\\
\noindent 93430 Villetaneuse, France\\
\vspace{10pt}

\noindent \textbf{Amie Wilkinson (wilkinso@math.northwestern.edu)}\\
\noindent Department of Mathematics, Northwestern University\\
\noindent 2033 Sheridan Road \\
\noindent Evanston, IL 60208-2730,  USA

\end{document}